\documentclass{article}

\usepackage{hyperref}

\usepackage{amsfonts}
\usepackage{amssymb}
\usepackage{textcomp}
\usepackage{latexsym,amsmath}
\usepackage{graphics}

\renewcommand{\phi}{\varphi}
\newcommand{\be}{\begin{equation}}
\newcommand{\ee}{\end{equation}}
\newcommand{\ba}{\begin{eqnarray}}
\newcommand{\ea}{\end{eqnarray}}
\newcommand{\ban}{\begin{eqnarray*}}
\newcommand{\ean}{\end{eqnarray*}}

\newcommand{\nul}{{\bf0}}
\newcommand{\rd}{{\mathbb R}^d}
\newcommand{\zd}{{\mathbb Z}^d}
\newcommand{\td}{{\mathbb T}^d}
\renewcommand{\r}{{\mathbb R}}
\newcommand{\z} {{\mathbb Z}}
\newcommand{\cn} {{\mathbb C}}
\newcommand{\n} {{\mathbb N}}

\newcommand{\ddd}{,\dots,}
\renewcommand{\lll}{\left(}
\newcommand{\rrr}{\right)}
\newcommand{\h}{\widehat}
\newcommand{\w}{\widetilde}
\newcommand{\too}{\mathop{\longrightarrow}}

\def\Z{{{\Bbb Z}}}

\def\R{{\Bbb R}}

\def\vp{{\varphi}}
\def\t{{\theta }}

\def\){\right)}
\def\({\left(}

\def\tt{\widetilde{t}}
\def\kt{\widetilde{k}}
\def\chit{\widetilde{\chi}}

\title{Differential and falsified sampling expansions
\thanks{This research was supported by Volkswagen Foundation; the first author is also
supported by H2020-MSCA-RISE-2014 Project number 645672; the second and the third authors are also
supported by grants from RFBR \# 15-01-05796-a,  St. Petersburg State University \#~9.38.198.2015.}}
\author{
Yu. Kolomoitsev$^{1, 2}$, A. Krivoshein$^{3}$ and M. Skopina$^{3}$
}
\date{ \small $^{1}$Universit\"at zu L\"ubeck,
Institut f\"ur Mathematik, L\"ubeck, Germany  \\
\small $^{2}$Institute of Applied Mathematics and Mechanics of NAS of Ukraine, Slov'yans'k, Ukraine
\\
 $^{3}$St. Petersburg State University, Russia \\
kolomoitsev@math.uni-luebeck.de, KrivosheinAV@gmail.com, skopina@ms1167.spb.edu}

\begin{document}
\maketitle

\begin{abstract}
Differential and falsified sampling expansions
$\sum_{k\in \mathbb{Z}^d}c_k\phi(M^jx+k)$, where $M$ is a matrix dilation, are studied.
In the case of differential expansions, $c_k=Lf(M^{-j}\cdot)(-k)$,
where $L$ is an appropriate differential operator.
For a large class of functions $\phi$, the approximation order of differential
expansions was recently studied.
Some smoothness of the Fourier transform of $\phi$ from this class is required.
In the present paper, we obtain similar results for a class of band-limited
functions $\phi$ with the discontinuous Fourier transform.
In the case of falsified expansions, $c_k$~is the mathematical expectation
of random integral average of a signal $f$ near the point $M^{-j}k$.
 To estimate the approximation order of the falsified
sampling expansions we compare them with the differential expansions.
Error estimations in $L_p$-norm are given in terms of the Fourier transform of $f$.
\end{abstract}




\newtheorem{theo}{Theorem}
\newtheorem{lem}[theo]{Lemma}
\newtheorem {prop} [theo] {Proposition}
\newtheorem {coro} [theo] {Corollary}
\newtheorem {defi} [theo] {Definition}
\newtheorem {rem} [theo] {Remark}
\newtheorem {ex} [theo] {Example}

\newcommand{\tocsecindent}{\hspace{0mm}}

\section{Introduction}

The well-known sampling theorem (Kotel'nikov's or Shannon's formula)
states that
\be
f(x)=\sum_{k\in\z}  f(-2^{-j}k)\, {\rm sinc}(2^jx+k), \quad{\rm sinc}(x):=\frac{\sin\pi x}{\pi x},
\label{0}
\ee
for  band-limited to $[-2^{j-1},2^{j-1}]$ signals (functions) $f$.
 This formula is very useful for engineers.
It was just Kotel'nikov~\cite{02} and Shannon~\cite{01} who started to
apply this formula for signal processing, respectively in 1933 and 1949.
Up to now, an overwhelming diversity of digital signal processing applications
and devices are based on it and  more than successfully use it.
 However, mathematicians knew this formula much earlier, actually, it can be found
in the papers by Ogura~\cite{03} (1920), Whittaker~\cite{W} (1915),
Borel~\cite{Borel} (1897), and even Cauchy~\cite{04} (1841).

Equality~(\ref0)  holds only for functions $f\in L_2(\r)$
whose Fourier transform is supported on $[-2^{j-1},2^{j-1}]$.
However the right hand side of~(\ref0) (the sampling expansion of $f$)
has meaning for every continuous $f$ with a good enough decay.
The problem of approximation of $f$ by its sampling expansions as $j\to+\infty$
was studied by many mathematicians. We  mention only some of such results.
 Brown~\cite{Brown} proved that  for every $x\in\r$
\be
\left|f(x)- \sum_{k\in\,\z} f(-2^{-j}k)\,{\rm sinc}(2^jx+k)\right|\le
C \int\limits_{|\xi|>2^{j-1}}
 |\h f(\xi)|d\xi,
\label{brown}
\ee
whenever the Fourier transform of  $f$ is summable on $\r$.
It is known that the pointwise approximation by sampling
expansions does not hold for arbitrary continuous functions $f$, even compactly supported.
Moreover, Trynin~\cite{Tr4} proved that there exists a continuous function vanishing
outside of $(0,\pi)$ such that its  deviation from the sampling expansion diverges
at every point $x\in(0,\pi)$.  Approximation by sampling expansions  in $L_p$-norm
was also  actively studied.
Bardaro,  Butzer,  Higgins, Stens, and Vinti~\cite{Butz4}, \cite{Butz5},
 proved that
$$
\Delta_p:=\Big\|f- \sum_{k\in\,\z} f(-2^{-j}k)\,{\rm sinc}(2^j\cdot+k)\Big\|_{L_p(\R)}
\too\limits_{j\to+\infty}0,\quad 1\le p<\infty,
$$
for $f\in C\cap \Lambda_p$, where $\Lambda_p$ consists of $f$
such that
$$
\sum_{k\in \Z}|f(x_k)|^p(x_{k}-x_{k-1})<\infty
$$
for some class of admissible
partition $\{x_k\}_{k\in \Z}$ of $\r$. Also they proved that the Sobolev spaces
${W_p^r}(\R)$, $r\in \n$, are subspaces of $\Lambda_p$, and that for every $f\in {W_p^r}(\R)$
\be
\Delta_p\le\frac{C\omega(f^{(r)}, 2^{-j})_p}{2^{-jr}},
\label{02}
\ee
where $\omega(\ )_p$ is the modulus of continuity in $L_p(\R)$, which is a typical
estimate for shift invariant spaces.
In~\cite{Butz6}, the same authors studied a generalized sampling approximation, replacing the
 sinc-function  by certain  linear combinations of $B$-splines.
For the case $p=\infty$,
Butzer,  Ries, and Stens~\cite{ButzNEW} proved that if a bounded function $\phi$ is such that
$$
\sup\limits_{x\in\r}\sum_{k\in \Z}|x-k|^{r+1}|\phi(x-k)|<\infty
$$
for some $r\in\n$, then the estimate
\be
\sup\limits_{x\in\r}\left|\sum_{k\in\Z} f\left(\frac kW\right)\phi(Wt-k)-f(t)\right|
\le C W^{-r}\omega\left(f^{(r)}, \frac 1W \right)
\label{02'}
\ee
holds for all $f\in C^r(\R)$ and $W>0$ if and only if the Strang-Fix condition of order
$r+1$  is satisfied for~$\phi$.
 The minimal requirement
for the convergence (in terms of Sobolev spaces) is $f\in W_p^1(\R)$, and, by (\ref{02}), the minimal
approximation order is $o(2^{-j})$.
 Sickel~\cite{Si1}, \cite{Si2}
studied error estimations for $\Delta_p$ in terms of Besov-Triebel-Lizorkin spaces. In particular,
his results provide the approximation order $2^{-j\gamma}$, $1/p\le\gamma<1$, for some functions.
The author of~\cite{Sk1}, \cite{Sk2} investigated approximation by sampling expansions
$$
\sum_{k\in\,\z} f(-2^{-j}k)\phi(2^jx+k)
$$
 for a wide class of band-limited functions $\phi$. For $p\ge2$,
 the error analysis was given in terms of the Fourier transform of $f$.
 In particular, the approximation order $2^{-j({1/p+\alpha})}$ was found for functions $f$
 in Sobolev spaces $W_1^1(\R)$ with $f'\in {\rm Lip}_{L_1} \alpha$, $\alpha>0$. In the case $1/p+\alpha>1$
these estimates give the same approximation order as~(\ref{02}).
Similar results were obtained for the generalized
sampling expansions (differential expansions)
$$
\sum_{k\in\,\z}  Lf(2^{-j}\cdot)(-k)\phi(2^jx+k),
$$
where $Lf:=\sum_{l=0}^m{\alpha_l}f^{(l)}$, and a function $\phi$ satisfies a special condition of compatibility with $L$.
The approximation order depends on $m$ in this case. An analog of
Brown's estimate~(\ref{brown}) was proved for such expansions in~\cite{Sk2}.

Multivariate  differential expansions
\be
\sum_{k\in\,\zd}  Lf(M^{-j}\cdot)(-k)\phi(M^jx+k),
\label{03}
\ee
where $M$ is a matric dilation, were studied in~\cite{KS1}.
Error estimations in terms of the Fourier transform of the approximated function were given
for a large class of functions $\phi$.
In particular, this class contains compactly supported functions, but it does not contain
functions with discontinuous Fourier transform.

Differential  expansions~(\ref{03})
 may be useful for  some problems and engineering applications.
The analog of Kotelnikov's formula~(\ref{0}) for differential  expansions
can be used to solve differential equations. The solution can be represented in
analytic form which depends only on sampled values of
known function for some equations, see Sec.~\ref{GenSamExp} for details.
On the other hand, the integral average value of a function $f$
near the point $M^{-j}k$ is close to $Lf(M^{-j}\cdot)(-k)$ for certain
operators $L$. Indeed, if $d=1$, $M=2$, then
\ban
\frac1{2^{-j}h}\int\limits_{2^{-j}k}^{2^{-j}(k+h)}f(t)\,dt\approx
\frac1{2^{-j}h}\int\limits_0^{2^{-j}h}\sum_{l=0}^N\frac1{l!}f^{(l)}(2^{-j}k)t^l\,dt=
\sum_{l=0}^N\frac1{(l+1)!}h^l\frac{d^lf(2^{-j}\cdot)}{dx^l}(k).
\ean
But the latter sum is nothing as $Lf(2^{-j}\cdot)(k)$, where
$\alpha_l=\frac1{(l+1)!}h^l$. This idea is used in Sec.~\ref{FalSamExp} for error
analysis of falsified sampling approximation.

In the present paper we study approximation properties of
differential  expansions~(\ref{03}) for a class of band-limited functions $\phi$
 with discontinuous Fourier transform.
Error estimations in $L_p$-norm are given in terms of the Fourier transform of the approximated function $f$.
Analogs of the classical sampling theorem and Brown's inequality~(\ref{brown})  are proved.
We also study the falsified sampling expansions
$$
\sum_{k\in\zd}E(f, M^{-j}k)\,\phi(M^jx+k),
$$
 where  $M$ is a matrix dilation,
$$
E(f, M^{-j}k)=E(f, M^{-j}k, h,w)=\int\limits_0^\infty du\,w(u)\frac  1 {{\rm meas }\,M^{-j}B_{h(u)}} \int\limits_{M^{-j}B_{h(u)}} f(M^{-j}k+t)d t,
$$
$h(u)$ is a positive function defined on $(0, \infty)$
and $u$ is a random value with probability density~$w$.
 To estimate the approximation order of falsified
sampling expansions we compare them with  differential expansions.
The error estimations of approximation by differential expansions
obtained in~\cite{KS1}  as well as new estimations for
band-limited functions $\phi$ are used.

The paper is organized as follows: in
Section~2 we introduce notation and give  some basic facts.
In Section~3 we study scaling operators $\sum_{k\in\zd} \langle f, {\w\phi}_{jk}\rangle \phi_{jk}$ and their approximation properties for a class of band-limited functions $\phi$
with discontinuous Fourier transform.
In Section~4 we study approximation properties of generalized sampling expansions defined by differential operators. In Section~5
we obtain estimates of the approximation order of falsified
sampling expansions. In Section~6 we give some examples.

\section{Notation and basic facts}
\label{notation}

$\n$ is the set of positive integers,
    $\r$  is the set of real numbers,
    $\cn$ is the set of complex numbers.
    $\rd$ is the
    $d$-dimensional Euclidean space,  $x = (x_1\ddd x_d)$ and $y =
    (y_1\ddd y_d)$ are its elements (vectors),
    $(x)_j=x_j$ for $j~=~1,\dots,d,$
    $(x, y)~=~x_1y_1+~\dots~+x_dy_d$,
    $|x| = \sqrt {(x, x)}$, ${\bf0}=(0\ddd 0)\in\rd$;
      $B_r=\{x\in\rd:\ |x|\le r\}$, $\td=[-\frac 12,\frac 12]^d$;
    $\zd$ is the integer lattice
    in $\rd$, $\z_+^d:=\{x\in\zd:~x\geq~{\bf0}\}.$
    If $\alpha,\beta\in\zd_+$, $a,b\in\rd$, we set
    $[\alpha]=\sum\limits_{j=1}^d \alpha_j$,
    $\alpha!=\prod\limits_{j=1}^d(\alpha_j!),$
    $$\binom{\beta}{\alpha}=\frac{\alpha!}{\beta!(\alpha-\beta)!},\quad
    a^b=\prod\limits_{j=1}^d a_j^{b_j},\quad
    D^{\alpha}f=\frac{\partial^{[\alpha]} f}{\partial x^{\alpha}}=\frac{\partial^{[\alpha]} f}{\partial^{\alpha_1}x_1\dots
    \partial^{\alpha_d}x_d},$$
    $\delta_{ab}$~is the Kronecker delta.

A real $d\times d$ matrix $M$ whose
eigenvalues are bigger than 1 in module is called a  dilation matrix.
Throughout the paper we
consider that such a matrix $M$ is fixed and  $m=|{\rm det}\,M|$,
$M^*$ denotes the conjugate matrix to $M$.
Since the spectrum of the operator $M^{-1}$ is
located in  $B_r$,
where $r=r(M^{-1}):=\lim_{j\to+\infty}\|M^{-j}\|^{1/j}$
is the spectral radius of $M^{-1}$, and there exists at least
one point of the spectrum on the boundary of $B_r$, we have
	\be
	\|M^{-j}\|\le {C_{M,\theta}}\,\vartheta^{-j},\quad j\ge0,
	\label{00}
	\ee
for every  positive number $\vartheta$  which is smaller in module
than any eigenvalue of  $M$.
In particular, we can take $\vartheta > 1$, then
	\be
	\lim_{j\to+\infty}\|M^{-j}\|=0.
	\label{000}
	\ee
A  matrix $M$ is called \emph{isotropic}
	if it is similar to a diagonal matrix
with elements $\lambda_1,\dots,\lambda_d$ that are placed on the main diagonal
	and $|\lambda_1|=\cdots=|\lambda_d|$.
	Thus, $\lambda_1,\dots,\lambda_d$ are eigenvalues of $M$
	and the spectral radius of $M$ is equal to $|\lambda|,$
	where $\lambda$ is one of the eigenvalues of $M.$
	Note that if the matrix $M$ is isotropic then
	$M^*$ is isotropic and $M^j$ is isotropic for all $j\in\z.$	

	It is well known that for an isotropic matrix $M$ and for  any $j\in\z$
		we have
	\be
	C_1 |\lambda|^j \le \|M^j\| \le C_2 |\lambda|^j,
	\label{10}
\ee
   where $\lambda$ is one of the eigenvalues of $M$ and the positive constants $C_1$ and  $C_2$ do not depend on $j$.

If $\phi$ is a function defined on $\rd$, we set
$$
\phi_{jk}(x):=m^{j/2}\phi(M^jx+k),\quad j\in\z,\,\, k\in\rd.
$$

$L_p$ denotes $L_p(\rd)$, $1\le p\le\infty$. We use $W_p^n$ (or $W_p^n(\R^d)$), $1\le p\le\infty$, $n\in\n$,
to denote  the Sobolev space on~$\rd$, i.e. the set of
functions whose derivatives up to order $n$ are in $L_p$,
with the usual Sobolev norm.

If $f, g$ are functions defined on $\rd$ and $f\overline g\in L_1$,
then  $\langle  f, g\rangle:=\int_{\rd}f\overline g$.

 For any function $f$, we set $f^{-}(x):=f(-x).$

If $F$ is a $1$-periodic (with respect to each variable) function and
$F\in L_1(\td)$,  then  $\widehat
F(k)=\int_{\td} F(x)e^{-2\pi i
(k,\,x)}\,dx$ is its $k$-th Fourier coefficient.
If $f\in L_1$,  then its Fourier transform is $\widehat
f(\xi)=\int_{\rd} f(x)e^{-2\pi i
(x,\xi)}\,dx$.

Denote by $\mathcal{S}$ the Schwartz class of functions defined on $\rd$.
    The dual space of $\mathcal{S}$ is $\mathcal{S}'$, i.e. $\mathcal{S}'$ is
    the space of tempered distributions.
    The basic facts from distribution theory
    can be found, e.g., in~\cite{Vladimirov-1}.
    Suppose $f\in \mathcal{S}$, $\phi \in \mathcal{S}'$, then
    $\langle \phi, f\rangle:= \overline{\langle f, \phi\rangle}:=\phi(f)$.
    If  $\phi\in \mathcal{S}',$  then $\h \phi$ denotes its  Fourier transform
    defined by $\langle \h f, \h \phi\rangle=\langle f, \phi\rangle$,
    $f\in \mathcal{S}$.
    If  $\phi\in \mathcal{S}'$, $j\in\z, k\in\zd$, then we define $\phi_{jk}$ by
        $
        \langle f, \phi_{jk}\rangle=
        \langle f_{-j,-M^{-j}k},\phi\rangle$ for all $f\in \mathcal{S}$.

Given a dilation  matrix $M$ and  $\delta>0$, we introduce a special notation for
the following integrals if they make sense
$$
{\cal I}_{j,\gamma,q}^{\rm In}(g):=
\int\limits_{|M^{*-j}\xi|<\delta}
	|\xi|^{q\gamma}|g(\xi)|^q d\xi, \quad
	{\cal I}_{j,\gamma,q}^{\rm Out}(g):=
\int\limits_{|M^{*-j}\xi|\ge\delta}
	|\xi|^{q\gamma}| g(\xi)|^q d\xi.
$$

A function $\phi\in L_1$ is said to satisfy the
{\em Strang-Fix condition of order $n$}
if  $D^{\beta} \h \phi (k)=0,$ whenever $k\in \zd\setminus\{\nul\}$
and $[\beta]<n.$

Let $1\le p \le \infty$. Denote by ${\cal L}_p$ the set
	$$
	{\cal L}_p:=
	\left\{
	\phi\in L_p\,:\, \|\phi\|_{{\cal L}_p}:=
	\left\|\sum_{k\in\zd} \left|\phi(\cdot+k)\right|\right\|_{L_p(\td)}<\infty
	\right\}.
	$$
	With the norm $\|\cdot\|_{{\cal L}_p}$, ${\cal L}_p$ is a Banach space.
	The simple properties are:
	${\cal L}_1=L_1,$
	$\|\phi\|_p\le \|\phi\|_{{\cal L}_p}$,
	$\|\phi\|_{{\cal L}_q}\le \|\phi\|_{{\cal L}_p}$
	for $1\le q \le p \le\infty.$ Therefore, ${\cal L}_p\subset L_p$
	and ${\cal L}_p\subset {\cal L}_q$ for $1\le q \le p \le\infty.$
	If $\phi\in L_p$ and compactly supported then $\phi\in {\cal L}_p$ for $p\ge1.$
	If $\phi$ decays fast enough, i.e. there exist constants $C>0$
	and $\varepsilon>0$ such that
	$|\phi(x)|\le C( 1+|x|)^{-d-\varepsilon}$ for all $x\in\rd,$
	then $\phi\in {\cal L}_\infty$.
	
	The following auxiliary statements will be useful for us.
	\begin{prop}[\cite{JiaMicPrewav}]	
\label{propLp}
	Let $1\le p \le \infty$. If $\phi \in {\cal L}_p$ and $a=\{a_k\}_{k\in\zd}\in \ell_p$, then
	$$
\left\|\sum_{k\in\zd} a_k \phi_{0k}\right\|_p\le
	\|\phi\|_{{\cal L}_p} \|a\|_{\ell_p}.
$$
\end{prop}

\begin{lem}[\cite{KS1}]
\label{lem1}
Let $1\le q <\infty$, $1/p+1/q=1$, $j\in\z_+$,
	$\w\phi$ be a tempered distribution
	whose Fourier transform $\h{\w\phi}$ is a function on $\rd$
	such that $|\h{\w\phi}(\xi)|\le C_{\w\phi} |\xi|^{N}$
	 for almost all $\xi\notin\td$, $N=N({\w\phi})\ge 0,$ and
	 $|\h{\w\phi}(\xi)|\le C'_{{\w\phi}}$
	 for almost all $\xi\in\td$.
	 Suppose $g\in L_q$, $g(\xi)=O(|\xi|^{-N-d-\varepsilon})$
as $|\xi|\to\infty$, where $\varepsilon>0$;  $\gamma\in
(N+\frac dp, N+\frac dp+\epsilon)$
 for $q\ne1$, $\gamma=N$ for $q=1$, and set
$$
	G_j(\xi)=G_j({\w\phi},g,\xi):=\sum\limits_{l\in\,\zd}
	g(M^{*j}(\xi+l))\overline{\h{\w\phi}(\xi+l)}.
$$
 	Then $G_j$ is a $1$-periodic function in
  	$L_q(\td)$, $\langle g,\h{{\w\phi}_{jk}}\rangle=m^{j/2}\h G_j(k)$,
  and for every $\delta\in(0,\frac 12)$
 	 \be
	\left\|G_j-g(M^{*j}\cdot)\overline{\h{\w\phi}}\right\|^q_{L_q(\td)}\le
    m^{- j} (C_{\gamma,\,{\w\phi}})^q \|M^{*-j}\|^{\gamma q} \,
   {\cal I}_{j,\gamma,q}^{\rm Out}(g).
	\label{fLem1GjLq}
 	\ee

\end{lem}

	\begin{lem}[\cite{KS1}]
\label{lemQjLp}
	Let $g$ and $\w\phi$  be as in Lemma~\ref{lem1}.
	Suppose $2\le p \le\infty$,
	 $\phi \in {\cal L}_p$.
	Then the series
	$\sum_{k\in\zd} \langle g,\h{\w\phi}_{jk}\rangle \phi_{jk}$
	converges unconditionally in $L_p$.
\end{lem}

\section{Scaling Approximation }
\label{scaleAppr}

Scaling operator
$\sum_{k\in\zd} \langle f, {\w\phi}_{jk}\rangle \phi_{jk}$
is a good tool of approximation for many  appropriate pairs of functions
$\phi, \w\phi$. We are interested in such operators, where $\w\phi$ is a
tempered distribution, e.g., the delta-function or a linear
combination of its derivatives. In this case the inner product
$\langle f, \w\phi_{jk}\rangle $ has meaning only for functions $f$ in $\mathcal{S}$.
To extend the class of functions $f$ one can  replace
$\langle f, {\w\phi}_{jk}\rangle $  by
$\langle \h f, \h{\w\phi_{jk}}\rangle$. In this case we set
$$
Q_j(\phi,\w\phi, f)=
 	\sum_{k\in\zd} \langle \widehat{f}, \h{\w\phi_{jk}}\rangle \phi_{jk},\quad j\in\z_+.
$$
Approximation properties of such operators for certain classes
of distributions~$\w\phi$ and functions $\phi$ were studied in~\cite{KS1}.
In particular, the following statement is proved.
	\begin{theo}[\cite{KS1}]
\label{theoQj}
	Let $2\le p \le \infty$, $1/p+1/q=1$, $N\in\z_+,$
$\gamma\in(N+\frac dp, N+\frac dp+\epsilon)$ for $p\ne\infty$,
 and $\gamma=N$ for $p=\infty$.
	Suppose
\begin{itemize}
  \setlength{\itemsep}{0cm}%
  \setlength{\parskip}{0cm}%

  \item[$(a)$]  $\w\phi$ is a tempered distribution
	 whose Fourier transform $\h{\w\phi}$ is a function on $\rd$
	such that $|\h{\w\phi}(\xi)|\le C_{\w\phi} |\xi|^{N}$
	 for almost all $\xi\notin\td$, $N>0,$ and
	 $|\h{\w\phi}(\xi)|\le C'_{{\w\phi}}$
	 for almost all $\xi\in\td$;

 	\item[$(b)$] $\phi \in {\cal L}_p $ and
	there exists $B_{\phi}>0$ such that
$
\sum_{k\in\zd}  |\h\phi(\xi+k)|^q<B_{\phi}$ for all $\xi\in\rd$;

	\item[$(c)$] there exist $n\in\n$ and $\delta\in(0, 1/2)$ such that
	$\h\phi\h{\w\phi}$ is  boundedly differentiable up to order $n$ on
	$\{|\xi|<\delta\}$,
	$\h\phi$ is boundedly	differentiable up to order $n$  on $\{|\xi+l|<\delta\}$	for all $l\in\zd\setminus\{\nul\}$;
	the function $\sum_{l\in\zd,\, l\neq\nul}|D^\beta \h \phi (\xi+l)|$
is bounded on	$\{|\xi|<\delta\}$  for $[\beta]=n$;
	$D^{\beta}(1-\h\phi\h{\w\phi})(0) = 0$ for $[\beta]<n$;
	the Strang-Fix condition of order $n$ holds 	for $\phi$;
	
      \item[$(d)$]  $f\in L_p$,
       $\h f\in L_q$, $\h f(\xi)=O(|\xi|^{-N-d-\varepsilon})$
as $|\xi|\to\infty$,
	$\varepsilon>0$.
 	
\end{itemize}

\noindent
	Then
	\be
	\|f-Q_j(\phi,\w\phi, f)\|_p^q\le
	C_1
	 \|M^{*-j}\|^{\gamma q}  {\cal I}_{j,\gamma,q}^{\rm Out}(\h f)+
	 C_2 \|M^{*-j}\|^{nq}   {\cal I}_{j,n,q}^{\rm In}(\h f),
	 \label{fTheoQjMain1}
	\ee
where $C_1$ and $C_2$ do not depend on $j$ and $f$.
\end{theo}	

A compability between  distribution $\w\phi$ and  function $\phi$
given in item (c) is required for the error of approximation in the latter theorem.
This compability is given in terms of the derivatives of $\h{\w\phi}$, $\h\phi$ at the origin
up to  the order of the Strang-Fix condition of $\phi$.
In what follows we will also consider a compability up to an arbitrary order.
\begin{defi}
\label{d1}
 A tempered distribution  $\w\phi$ and a function $\phi$ is said to be
 {\em strictly compatible} if there exists $\delta\in(0,1/2)$ such that
 $\overline{\h\phi}(\xi)\h{\w\phi}(\xi)=1$
 a.e. on $\{|\xi|<\delta\}$ and $\h\phi(\xi)=0$ a.e.
 on $\{|l-\xi|<\delta\}$ for all
 $l\in\z\setminus\{0\}$.
\end{defi}

The class of functions $\phi$ considered in Theorem~\ref{theoQj} is large,
it contains both compactly supported and band-limited functions.
However, it does not include functions whose Fourier
transform is discontinuous, for example, the function $\prod_{k=1}^d{\rm sinc}(x_k)$
is out of consideration. Here we will make up  for this omission.

 Let ${\cal B}={\cal B}(\R^d)$  denote the class of functions $\phi$ given by
 \be
\phi(x)=\int\limits_{\rd}\theta(\xi)e^{2\pi i(x,\xi)}\,d\xi,
\label{61}
\ee
where $\theta$ is supported on a parallelepiped $S:=[a_1, b_1]\times\dots\times[a_d, b_d]$ and such that	$\theta\big|_S\in C^d(S)$.

\begin{prop}
\label{prop1}
	Let $1< p < \infty$, $\phi\in \cal B$, $f\in L_p$. Then
	\begin{equation}\label{eqK1}
  \left(\sum_{k\in\zd} |\langle f,{{\phi}_{0k}}\rangle|^p\right)^\frac 1p\le C_{\phi, p}\|f\|_p.
\end{equation}
\end{prop}
 The proof of Proposition~\ref{prop1} is based on  two lemmas.
To formulate these lemmas, as well as further results, we need additional notation. Set
$$
U_k^0=\{t\in \r\,:\,|t-k|<1\}\quad\text{and}\quad U_k^1=\r\setminus U_k^0,\quad k\in \z;
$$
if $k\in\zd$, $\chi=(\chi_1,\dots,\chi_d)\in \{0,1\}^d$, then $U_k^\chi$ is defined by
$$
U_k^\chi = U_{k_1}^{\chi_1}\times\dots\times U_{k_d}^{\chi_d}.
$$

The proof of the lemma below follows easily from the proof of Lemma~4 in~\cite{Sk1}.

\begin{lem}\label{lemK1}
  Let $f\in L_p(\r)$, $1<p<\infty$, $u\in \r$. Then
  $$
   \lll\sum_{k\in \z}\bigg|\,\,\int\limits_{U_k^1}f(t)\frac{e^{2\pi i u(t-k)}}{t-k}dt\bigg|^p\rrr^\frac1p \le C_p \|f\|_{L_p(\r)}.
  $$
\end{lem}

The following lemma gives a  proof of Proposition~\ref{prop1} in the case $d=1$. Since $S=[a_1,b_1]$ if $d=1$, for convenience,
we redenote $[a_1, b_1]$ by $[a, b]$ for this case.

\begin{lem}\label{lemK2}
  Let $f\in L_p(\r)$, $1<p<\infty$, and $\phi\in \mathcal{B}(\R)$. Then, for any $\chi\in \{0,1\}$, we have
\begin{equation}\label{eqK2}
  \lll\sum_{k\in \z} \bigg|\int\limits_{U_k^\chi} f(t) \phi(t-k) dt\bigg|^p \rrr^\frac1p \le C_{S,p}\Vert \t\Vert_{W_\infty^1(S)}\|f\|_{L_p(\r)}.
\end{equation}
\end{lem}
{\bf Proof. }
1) If $\chi=0$, then by H\"older's inequality, we have
\begin{equation*}\label{eqK3}
  \begin{split}
     \sum_{k\in \z} \bigg|\int\limits_{U_k^0} f(t) \phi(t-k) dt\bigg|^p \le &\|\phi\|_{L_\infty(S)}^p \sum_{k\in\z}\lll \, \int\limits_{U_k^0}|f(t)|dt\rrr^p \le\\
     & C_{S,p}\Vert\t\Vert_{L_\infty(S)}^p \sum_{k\in\z}\,\int\limits_{U_k^0}|f(t)|^p dt\le C_{S,p}\Vert\t\Vert_{L_\infty(S)}^p\|f\|_{L_p(\r)}^p,
  \end{split}
\end{equation*}
which implies~\eqref{eqK2}.

2) Let $\chi=1$. Using the formula
\begin{equation}\label{eqK4}
  \phi(x)=\int_a^b \theta(\xi)e^{2\pi i x \xi}d\xi=\frac{\theta(b)e^{2\pi ibx}-\theta(a)e^{2\pi iax}}{2\pi ix}-\frac{1}{2\pi ix}\int\limits_a^b \theta'(\xi)e^{2\pi x\xi}d\xi,
\end{equation}
we obtain
\begin{equation}\label{eqK5}
  \begin{split}
    \sum_{k\in\Z}\bigg| \int\limits_{U_k^1} f(t)\vp(t-k)dt\bigg|^p\le &C_p\sum_{k\in\Z}\bigg| \int\limits_{U_k^1}f(t)\frac{\theta(b)e^{2\pi ibx}-\theta(a)e^{2\pi iax}}{t-k}dt\bigg|^p+\\
    &C_p\sum_{k\in\Z}\bigg|\int\limits_{U_k^1}\frac{f(t)}{t-k}\bigg(\int_a^b \theta'(\xi)e^{2\pi i(t-k)\xi}d\xi\bigg)dt\bigg|^p=I_1+I_2.
  \end{split}
\end{equation}
By Lemma~\ref{lemK1}, we get
\begin{equation}\label{eqK6}
  I_1\le C_p\sum_{u\in\{a,b\}}|\theta(u)|^p \sum_{k\in\Z}\bigg|\int\limits_{U_k^1}f(t)\frac{e^{2\pi iu(t-k)}}{t-k} dt\bigg|^p\le C_{p}\Vert\t\Vert_{L_\infty(S)}^p \Vert f\Vert_{L_p(\R)}^p.
\end{equation}
Now, let us consider the sum $I_2$. Using H\"older's inequality with $1/p+1/q=1$ and Lemma~\ref{lemK1}, we derive
\begin{equation}\label{eqK7}
  \begin{split}
    I_2=&C_p\sum_{k\in\Z}\bigg|\int\limits_a^b \theta'(\xi)\int\limits_{U_k^1}f(t)\frac{e^{2\pi i(t-k)\xi}}{t-k} dtd\xi\ \bigg|^p\le\\
    &C_p \sum_{k\in\Z}\(\int\limits_a^b |\theta'(\xi)|\bigg|\int\limits_{U_k^1}f(t)\frac{e^{2\pi i(t-k)\xi}}{t-k} dt\bigg|d\xi\ \)^p\le\\
    &C_p\sum_{k\in\Z}\(\int\limits_a^b|\theta'(\xi)|^q d\xi\)^\frac pq \int\limits_a^b \bigg|\int\limits_{U_k^1}f(t)\frac{e^{2\pi i(t-k)\xi}}{t-k} dt\bigg|^p d\xi\le\\
    &C_{S,p}\Vert \t'\Vert_{L_\infty(S)}^p\int\limits_a^b \sum_{k\in\Z}\bigg|\int\limits_{U_k^1}f(t)\frac{e^{2\pi i(t-k)\xi}}{t-k} dt\bigg|^p d\xi\le C_{S,p}\Vert \t'\Vert_{L_\infty(S)}^p \Vert f\Vert_{L_p(\R)}^p.
  \end{split}
\end{equation}

Finally, combining \eqref{eqK5}--\eqref{eqK7}, we get~\eqref{eqK2} for $\chi=1$.
$\diamond$

\bigskip

{\bf Proof of Proposition~\ref{prop1}.} For convenience, we introduce the following notation.
If $t\in\rd$, then  $\widetilde{t}:=(t_1,\dots,t_{d-1})\in \R^{d-1}$.
For a function $g$ of $d$ variables $t_1,\dots t_{d-1}$ and $s$, we set $g_s(\widetilde{t})=g(t_1,\dots,t_{d-1},s)$.
 Let also $\psi_{\widetilde{x}}(\eta)=\mathcal{F}^{-1}\theta_\eta(\widetilde{x})$ and $\widetilde{S}=[a_1, b_1]\times\dots\times[a_{d-1}, b_{d-1}]$.

We have
\begin{equation}\label{eqK8}
  \begin{split}
    \sum_{k\in\Z^d}|\langle f,\vp_{0k} \rangle|^p=\sum_{k\in\Z}\bigg|\int\limits_{\R^d} f(t)\vp(t-k)dt \bigg|^p\le C_p \sum_{\chi\in\{0,1\}^d}I^\chi,
  \end{split}
\end{equation}
where
$$
I^\chi=\sum_{k\in\Z}\bigg|\int\limits_{U_k^\chi} f(t)\vp(t-k)dt \bigg|^p.
$$

Thus, to prove~\eqref{eqK1} it is enough to show that
\begin{equation}\label{eqK9}
 I^\chi \le C_{S,p}\Vert \t\Vert_{W_\infty^d(S)}^p\Vert f\Vert_p^p\quad\text{for any}\quad \chi\in \{0,1\}^d.
\end{equation}

We  prove~\eqref{eqK9} by  induction on $d$. For $d=1$, the inequality~\eqref{eqK9} was proved in Lemma~\ref{lemK2}.
To prove the inductive step $d-1\to d$, we assume that for any $g\in L_p(\R^{d-1})$ and $\widetilde{\phi}\in \mathcal{B}(\R^{d-1})$
(more precisely $\widetilde{\phi}=\mathcal{F}^{-1}\widetilde{\t}$, where $\widetilde{\t}$ is the same as in~\eqref{61} with $\widetilde{S}$ in place of $S$), we have
\begin{equation}\label{eqK10}
  \begin{split}
    \sum_{\widetilde{k}\in\Z^{d-1}}  \bigg|  \int\limits_{U_{\kt}^{\chit}} g\(\tt\)\widetilde{\vp}\(\tt-\kt\)d\tt  \bigg|^p\le C_{\widetilde{S},p}\Vert \widetilde{\t}\Vert_{W_\infty^{d-1}(\widetilde{S})}^p\Vert g\Vert_{L_p(\R^{d-1})}^p.
  \end{split}
\end{equation}

For any $\chi\in \{0,1\}^d$, we can write $\chi=(\chit,\chi_d)$ and
$$
I^\chi=\sum_{\chi_d\in \{0,1\}}I^{(\chit,\chi_d)}.
$$
Let us estimate $I^{(\chit,\chi_d)}$ for $\chi_d=0$ and $\chi_d=1$.

1) For $\chi_d=0$, we obtain
\begin{equation}\label{eqK12}
  I^{(\chit,0)}=\sum_{\kt\in\Z^{d-1}}\sum_{l\in\Z}\bigg| \int\limits_{U_l^0}\int\limits_{U_{\kt}^{\chit}} f_s\(\tt\){\vp}\(\tt-\kt,s-l\)d\tt ds\bigg|^p.
\end{equation}
Using the above notation and formula~\eqref{eqK4}, we get
\begin{equation}\label{eqK7.5}
  \begin{split}
    \vp(x)=&\mathcal{F}^{-1}\t(x)=\int\limits_{a_d}^{b_d} \psi_{\widetilde{x}}(\eta)e^{2\pi i x_d\eta}d\eta=\\
    &\frac{\psi_{\widetilde{x}}(b_d)e^{2\pi ib_dx_d}-\psi_{\widetilde{x}}(a_d)e^{2\pi ia_dx_d}}{2\pi ix_d}-\frac1{2\pi ix_d}\mathcal{F}^{-1}\psi'_{\widetilde{x}}(x_d).
  \end{split}
\end{equation}
By~\eqref{eqK7.5} and H\"older's inequality, we derive
\begin{equation}\label{eqK13}
  \begin{split}
    &\bigg|\int\limits_{U_l^0}\int\limits_{U_{\kt}^{\chit}} f_s(\tt)\int_{a_d}^{b_d} \psi_{\tt-\kt}(\eta)e^{2\pi i(s-l)\eta} d\eta d\tt ds\bigg|^p\le \(\int\limits_{a_d}^{b_d}\int\limits_{U_l^0}\bigg|\int\limits_{U_{\kt}^{\chit}}f_s(\tt)\psi_{\tt-\kt}(\eta) d\tt  \, \bigg|dsd\eta\)^p\le\\
    &\(\int\limits_{a_d}^{b_d}\int\limits_{U_l^0} dsd\eta\)^{p-1} \int\limits_{a_d}^{b_d}\int\limits_{U_l^0}\bigg|\int\limits_{U_{\kt}^{\chit}}f_s(\tt)\psi_{\tt-\kt}(\eta) d\tt  \,\bigg|^p dsd\eta
    =(2(b_d-a_d))^{p-1} \int\limits_{a_d}^{b_d}\int\limits_{U_l^0}|F_{\kt,\eta}(s)|^p dsd\eta,
  \end{split}
\end{equation}
where
\begin{equation}\label{eqKFun}
  F_{\kt,\eta}(s)=\int\limits_{U_{\kt}^{\chit}} f_s(\tt)\psi_{\tt-\kt}(\eta) d\tt.
\end{equation}

Next, combining~\eqref{eqK12} and~\eqref{eqK13} and using the induction hypothesis~\eqref{eqK10}, we get
\begin{equation}\label{eqK14}
  \begin{split}
    I^{(\chit,0)}\le &C_{S,p}\sum_{\kt\in\Z^{d-1}}\sum_{l\in\Z} \int\limits_{a_d}^{b_d}\int\limits_{U_l^0}|F_{\kt,\eta}(s)|^p dsd\eta
    \le C_{S,p}\int\limits_{a_d}^{b_d}\sum_{\kt\in\Z^{d-1}} \|F_{\kt,\eta}\|_{L_p(\R)}^p d\eta\le \\
    &C_{S,p}\int\limits_{a_d}^{b_d}\int\limits_\R \sum_{\kt\in\Z^{d-1}}\bigg| \int_{U_{\kt}^{\chit}} f_s\(\tt\)\mathcal{F}^{-1}\t_\eta\(\tt-\kt\)d\tt\bigg|^pdsd\eta\le \\
    &C_{S,p}\int\limits_{a_d}^{b_d} \Vert \t_\eta\Vert_{W_\infty^{d-1}(\widetilde{S})}^p d\eta\int\limits_\R \Vert f_s\Vert_{L_p(\R^{d-1})}^p ds\le C_{S,p} \Vert \t\Vert_{W_\infty^{d}(S)}^p\Vert f\Vert_p^p.
  \end{split}
\end{equation}

2) Now we consider the case $\chi_d=1$. We have
\begin{equation}\label{eqK17}
  I^{(\chit,1)}=\sum_{\kt\in\Z^{d-1}}\sum_{l\in\Z}\bigg| \int\limits_{U_l^1}\int\limits_{U_{\kt}^{\chit}} f_s\(\tt\){\vp}\(\tt-\kt,s-l\)d\tt ds\bigg|^p.
\end{equation}
By~\eqref{eqK7.5}, we get
\begin{equation}\label{eqK18}
  I^{(\chit,1)}\le C_p(I_1+I_2),
\end{equation}
where
\begin{equation*}
  I_1=\sum_{\kt\in\Z^{d-1}}\sum_{l\in\Z}\bigg| \int\limits_{U_l^1}\int\limits_{U_{\kt}^{\chit}} f_s\(\tt\)\frac{\psi_{\tt-\kt}(b_d)e^{2\pi ib_d(s-l)}-\psi_{\tt-\kt}(a_d)e^{2\pi ia_d(s-l)}}{s-l}d\tt ds \bigg|^p
\end{equation*}
and
\begin{equation}\label{eqK19}
  I_2=\sum_{\kt\in\Z^{d-1}}\sum_{l\in\Z}\bigg| \int\limits_{U_l^1}\int\limits_{U_{\kt}^{\chit}} f_s\(\tt\) \frac{\mathcal{F}^{-1}\psi'_{\tt-\kt}(s-l)}{s-l}d\tt ds \bigg|^p.
\end{equation}

Using the function~\eqref{eqKFun} and Lemma~\ref{lemK1}, we obtain
\begin{equation}\label{eqK19.5}
  \begin{split}
     I_1\le C_p\sum_{u\in\{a_d,b_d\}}\sum_{\kt\in\Z^{d-1}}\sum_{l\in\Z} \bigg|\int\limits_{U_l^1}\frac{F_{\kt,u}(s)e^{2\pi iu(s-l)}}{s-l}ds \bigg|^p\le  C_p\sum_{u\in\{a_d,b_d\}}\sum_{\kt\in\Z^{d-1}} \Vert F_{\kt,u}\Vert_{L_p(\R)}^p.
  \end{split}
\end{equation}
Now, using the induction hypothesis~\eqref{eqK10}, for any $u\in \R$, we derive
\begin{equation}\label{eqK19.6}
  \begin{split}
     \sum_{\kt\in\Z^{d-1}} \Vert F_{\kt,u}\Vert_{L_p(\R)}^p=&\int\limits_\R\sum_{\kt\in\Z^{d-1}} \bigg|\int\limits_{U_{\kt}^{\chit}} f_s(\tt)\mathcal{F}^{-1}\t_u(\tt-\kt) d\tt \bigg|^p ds\le\\
     &C_{\widetilde{S},p}\Vert \t_u\Vert_{W_\infty^{d-1}(\widetilde{S})}^p \int\limits_\R \Vert f_s\Vert_{L_p(\R^{d-1})}^p ds.
  \end{split}
\end{equation}
Hence, combining~\eqref{eqK19.5} and~\eqref{eqK19.6}, we get
\begin{equation}\label{eqK20}
  I_1\le C_{\widetilde{S},p} \sum_{u\in\{a_d,b_d\}}\Vert \t_u\Vert_{W_\infty^{d-1}(\widetilde{S})}^p \Vert f\Vert_p^p \le C_{S,p}\Vert \t\Vert_{W_\infty^{d}(S)}^p\Vert f\Vert_p^p.
\end{equation}

Let us consider $I_2$. Denoting
$$
F_{\kt,\eta}^*(s)=\int\limits_{U_{\kt}^{\chit}} f_s(\tt)\psi_{\tt-\kt}'(\eta) d\tt
$$
and using H\"older's inequality, we obtain
\begin{equation}\label{eqK21}
  \begin{split}
    \bigg|\int\limits_{U_l^1}\int\limits_{U_{\kt}^{\chit}} f_s(\tt)\frac{\mathcal{F}^{-1}\psi'_{\tt-\kt}(s-l)}{s-l}d\tt ds  \bigg|^p
    \le &\(\int\limits_{a_d}^{b_d}\bigg|\int\limits_{U_l^1}\int\limits_{U_{\kt}^{\chit}} f_s(\tt) \frac{\psi'_{\tt-\kt}(\eta)e^{2\pi i(s-l)\eta}}{s-l}d\tt ds\bigg| d\eta\)^p\le\\
    &(b_d-a_d)^{p-1} \int\limits_{a_d}^{b_d}\bigg|\int\limits_{U_l^1}F_{\kt,\eta}^*(s)\frac{e^{2\pi i(s-l)\eta}}{s-l} ds \bigg|^pd\eta.
  \end{split}
\end{equation}
Thus, combining~\eqref{eqK19} and~\eqref{eqK21}, using Lemma~\ref{lemK1}, and the induction hypothesis~\eqref{eqK10}, we derive
\begin{equation}\label{eqK22}
  \begin{split}
    I_2 &\le C_{S,p}\sum_{\kt\in\Z^{d-1}}\sum_{l\in\Z} \int\limits_{a_d}^{b_d}\bigg|\int\limits_{U_l^1}F_{\kt,\eta}^*(s)\frac{e^{2\pi i(s-l)\eta}}{s-l} ds \bigg|^pd\eta=\\
    &C_{S,p}\int\limits_{a_d}^{b_d}\sum_{\kt\in\Z^{d-1}}\sum_{l\in\Z} \bigg|\int\limits_{U_l^1}F_{\kt,\eta}^*(s)\frac{e^{2\pi i(s-l)\eta}}{s-l} ds \bigg|^pd\eta\le\\
    &C_{S,p}\int\limits_{a_d}^{b_d}\sum_{\kt\in\Z^{d-1}} \Vert F_{\kt,\eta}^* \Vert_{L_p(\R)}^p d\eta=\\
    &C_{S,p}  \int\limits_{a_d}^{b_d}\int\limits_\R \sum_{\kt\in\Z^{d-1}} \bigg|\int\limits_{U_{\kt}^{\chit}} f_s(\tt)\mathcal{F}^{-1}\frac{\partial}{\partial \eta}\t_\eta\(\tt-\kt\)d\tt \bigg|^p dsd\eta\le\\
    &C_{S,p} \int\limits_{a_d}^{b_d}\bigg\Vert \frac{\partial}{\partial \eta}\t_\eta\bigg\Vert_{W_\infty^{d-1}(\widetilde{S})}^pd\eta
    \int\limits_\R \Vert f_s \Vert_{L_p(\R^{d-1})}^p ds\le C_{S,p} \Vert \t\Vert_{W_\infty^{d}(S)}^p \Vert f\Vert_p^p.
  \end{split}
\end{equation}

Finally, combining~\eqref{eqK18}, \eqref{eqK20}, and~\eqref{eqK22}, we get~\eqref{eqK9} for any $\chi=(\chi_1,\dots,\chi_{d-1},1)$. This and~\eqref{eqK14} prove the proposition.
$\diamond$

\bigskip

\begin{prop}
\label{prop2}
	Let $1< p< \infty$,
	$\phi\in \cal B$, $a=\{a_k\}_{k\in\zd}\in\ell_p$. Then
	$$
\left\|\sum_{k\in\zd} a_k \phi_{0k}\right\|_p\le C_{\phi, p}\|a\|_{\ell_p}.
$$
\end{prop}

{\bf Proof. } Since
$$
\left\|\sum_{k\in\zd} a_k \phi_{0k}\right\|_p=\left|\left\langle\sum_{k\in\zd} a_k \phi_{0k}, f\right\rangle\right|=
\left|\sum_{k\in\zd} a_k \langle\phi_{0k}, f\rangle\right|,
$$
where $f\in L_q$, $\|f\|_q\le1$, $1/p+1/q=1$,  the statement follows immediately from Proposition~\ref{prop1} and H\"older's inequality. $\diamond$

\begin{coro}
\label{coro_prop2}
	Let $2\le p < \infty$, $1/{p}+1/{q}=1$, $g\in L_q$, and $\w\phi$  be as in Lemma~\ref{lem1},
	$\phi\in \cal B$. 		Then the series
	$\sum_{k\in\zd} \langle g,\h{\w\phi_{jk}}\rangle \phi_{jk}$
	converges unconditionally in $L_p$ and
	\be
	\left\|\sum_{k\in\zd} \langle g,\h{\w\phi_{jk}}\rangle \phi_{jk}\right\|_p\le C_{\phi, q}
	m^{\frac j2-\frac jp}\left(\sum_{k\in\zd}
	|\langle g,\h{{\w\phi}_{jk}}\rangle|^p\right)^{\frac 1p}.
	\label{0_prop}
	\ee
\end{coro}

{\bf Proof. } Because of Lemma~\ref{lem1} and the  Hausdorff-Young inequality, we
have
	\be
	\left(\sum_{k\in\zd} |\langle g,\h{{\w\phi}_{jk}}\rangle|^p\right)^\frac 1p =
	m^{\frac j2}
	\left(\sum_{k\in\zd} |\h G_j(k)|^p\right)^\frac 1p \le
	m^{\frac j2}  \|G_j\|_{L_q(\td)} <\infty,
	\label{1_prop}
	\ee
where $G_j$ is a function from Lemma~\ref{lem1}.
By Proposition~\ref{prop2},  we can state that for every
finite	subset $\Omega$ of $\zd$
	$$\left\|\sum_{k\in\Omega}
	 \langle g,\h{{\w\phi}_{jk}}\rangle \phi_{jk}\right\|_p=
m^{\frac j2-\frac jp}\left\|\sum_{k\in\Omega}
	 \langle g,\h{{\w\phi}_{jk}}\rangle \phi_{0k}\right\|_p
	\le
	m^{\frac j2-\frac jp}C_{\phi, q}
	\left(\sum_{k\in\Omega}
	|\langle g,\h{{\w\phi}_{jk}}\rangle|^p\right)^{\frac 1p}.
$$
	The series $\sum_{k\in\zd} |\langle g,\h{{\w\phi}_{jk}}\rangle|^p$ is
	convergent, which yields that
	$\sum_{k\in\zd} \langle g,\h{{\w\phi}_{jk}}\rangle \phi_{jk}$
	converges unconditionally. Similarly we obtain~(\ref{0_prop})	$\diamond$

\bigskip

Now we are ready to study approximation properties of the operators $$Q_j(\phi,\w\phi, f)=\sum_{k\in\zd} \langle \widehat{f}, \h{\w\phi_{jk}}\rangle \phi_{jk},$$ where $\phi\in \cal B$.
The following theorem is a counterpart of Theorem~\ref{theoQj} for such functions $\phi$.

\begin{theo}
\label{theoQj_new}
	Let $2\le p < \infty$, $1/p+1/q=1$, $N\in\z_+,$
$\gamma\in(N+\frac dp, N+\frac dp+\epsilon)$.
	Suppose
\begin{itemize}
  \setlength{\itemsep}{0cm}%
  \setlength{\parskip}{0cm}%

\item[$(a)$]  $\w\phi$ is a tempered distribution
	 whose Fourier transform $\h{\w\phi}$ is a function on $\rd$
	such that $|\h{\w\phi}(\xi)|\le C_{\w\phi} |\xi|^{N}$
	 for almost all $\xi\notin\td$ and
	 $|\h{\w\phi}(\xi)|\le C'_{{\w\phi}}$
	 for almost all $\xi\in\td$;

\item[$(b)$] $\phi\in \cal B$;
	
\item[$(c)$] $\w\phi$ and $\phi$ are strictly compatible;

\item[$(d)$]  $f\in L_p$,
       $\h f\in L_q$, $\h f(\xi)=O(|\xi|^{-N-d-\varepsilon})$
as $|\xi|\to\infty$,  $\varepsilon>0$.

\end{itemize}

\noindent
	Then
	\be
	\|f-Q_j(\phi,\w\phi, f)\|_p^q\le C
	 \|M^{*-j}\|^{\gamma q}  {\cal I}_{j,\gamma,q}^{\rm Out}(\h f),
	 \label{8}
	\ee
where $C$ does not depend on $j$ and $f$.
\end{theo}

{\bf Proof.}
Throughout the proof we denote by $C_1, C_2,\dots$
different constants which do not depend on $j$ and $f$.
It follows from~(\ref{0_prop}), (\ref{1_prop}),
and Lemma~\ref{lem1} with $g=\widehat{f}$  that
	\ban
	\|Q_j(\phi,\w\phi, f)\|_p\le	m^{\frac jq}  C_{\phi,q}\|G_j\|_{L_q(\td)}\le	
	C_{\phi,q} m^{\frac jq}\lll m^{- j} (C_{\gamma,\,{\w\phi}})^q \|M^{*-j}\|^{\gamma q} {\cal I}_{j,\gamma,q}^{\rm Out}(\h f)\rrr^{1/q}+
	\\
	 C_{\phi,q}m^{\frac jq}\Big(\int\limits_{\td}|\h f(M^{*j}\xi)\h{\w\phi}(\xi)|^q\,d\xi\Big)^{1/q}\le
	C_1\lll\Big(\|M^{*-j}\|^{\gamma q} {\cal I}_{j,\gamma,q}^{\rm Out}(\h f)\Big)^{1/q}+\|\h f\|_q\rrr.
	\ean
Since $f^-$ coincides with  $\h {\h f}$ a.e. due to the du Bois-Reymond lemma, applying the Hausdorff-Young inequality, we have
     $      \|f\|_p= \|f^{-}\|_p\le      \|\h f\|_q$, which yields
\ba
	\|f-Q_j(\phi,\w\phi, f)\|_p\le C_2\lll\Big(\|M^{*-j}\|^{\gamma q} {\cal I}_{j,\gamma,q}^{\rm Out}(\h f)\Big)^{1/q}+
	\|\h f\|_q\rrr
	\le
	\nonumber
	\\
	\nonumber
	C_2 \lll\|M^{*-j}\|^{\gamma q} {\cal I}_{j,\gamma,q}^{\rm Out}(\h f)+\int\limits_{|M^{*-j}\xi|<\delta}|
\h f(\xi)|^qd\xi+\delta^{-q\gamma}\int\limits_{|M^{*-j}\xi|
\ge\delta}|M^{*-j}\xi|^{q\gamma}|\h f(\xi)|^qd\xi\rrr^{1/q}\le
	\\
	C_3 \lll\|M^{*-j}\|^{\gamma q} {\cal I}_{j,\gamma,q}^{\rm Out}(\h f)+\int\limits_{|M^{*-j}\xi|<\delta}|
\h f(\xi)|^qd\xi\rrr^{1/q},
\label{35}
\ea		
where $\delta$ is from Definition~\ref{d1}.
		
		For any compact set $K\subset \rd$, the function $\h f$ can be approximated
 in $L_q(K)$ by infinitely smooth functions  supported on $K$.
So, given $j$, one can find a function $F_j$ such that
${\rm supp}\,\h F_j\subset \{|M^{*-j}\xi|<\delta\}$, $\h F_j\in C^\infty(\rd)$, and
\ban
\int\limits_{|M^{*-j}\xi|<\delta}|
\h f(\xi)-\h F_j(\xi)|^qd\xi\le
 \|M^{*-j}\|^{\gamma q}  {\cal I}_{j,\gamma,q}^{\rm Out}(\h f).
\ean
Combining this with~(\ref{35}), where
 $f$ is replaced by $f-F_j$, and taking into account that
${\cal I}_{j,\gamma,q}^{\rm Out}(\h f)={\cal I}_{j,\gamma,q}^{\rm Out}(\h f-\h F_j)$, we have
\ba
\|(f-F_j)-Q_j(\phi,\w\phi, f-F_j)\|_p^q\le
\|M^{*-j}\|^{\gamma q}  {\cal I}_{j,\gamma,q}^{\rm Out}(\h f).
\label{36}
\ea
Due to  Carleson's theorem
and Lemma~\ref{lem1},  we have
$$
\sum\limits_{k\in\zd}\langle \h F_j,\h{{\w\phi}_{jk}}
\rangle\widehat{\phi_{jk}}(\xi)=
\sum\limits_{l\in\,\zd}
\h F_j(\xi+M^{*j}l)\overline{\h{\w\phi}(M^{*-j}\xi+l)}\h{\phi}(M^{*-j}\xi).
$$
If $l\ne0$ and $F_j(\xi+M^{*j}l)\ne0$, then $|M^{*-j}\xi+l|<\delta$ and hence
$\h{\phi}(M^{*-j}\xi)=0$. So,
$$
\sum\limits_{k\in\zd}\langle \h F_j,\h{{\w\phi}_{jk}}
\rangle\widehat{\phi_{jk}}(\xi)=\h F_j(\xi),
$$
which yields that $Q_j(\phi,\w\phi, F_j)=F_j$.
It follows that
$$
\|f- Q_j(\phi,\w\phi, f)\|_p=\|f-F_j- Q_j(\phi,\w\phi, f-F_j)\|_p.
$$
Together with~(\ref{36}) this yields~(\ref{8}).
 $\Diamond$	

\begin{coro}
If, under  assumptions $(a)-(c)$ of Theorem~\ref{theoQj_new}, a function $f$ is such that
its Fourier transform is supported in $\{\xi:\ |M^{*-j}\xi|<\delta\}$,
where $\delta$ is from Definition~\ref{d1}, then
$$
f=Q_j(\phi,\w\phi, f)\quad a.e.
$$
\end{coro}	

The latter equality is a generalization of~(\ref{0}).	
	
\section{Differential expansions}
\label{GenSamExp}	

Consider  a differential operator $L$
	defined by
	\be
Lf:=\sum_{[\beta]\le N} a_\beta D^{\beta} f,
\quad a_{\beta}\in{\mathbb C},\,\, a_{\nul}\ne0,
\label{11}
\ee
where	$N\in\z_+.$
	The action of the operator $D^{\beta}$ is associated with the action of the corresponding derivative of the $\delta$-function.
	In  more detail, let $f$ be such that $
	\int_{\rd} (1+|\xi|)^{N+\alpha} |\h f(\xi)| d\xi < \infty,
	$
	$\alpha>0$,
	which implies that  $f$ is  continuously differentiable on $\rd$
	up to the order $N.$
	Then for  $[\beta]\le N $
	\ban
	D^{\beta}f (M^{-j}\cdot)(-k)=
	(-1)^{[\beta]} m^j \int\limits_{\rd} \h f(M^{*j}\xi) (-2\pi i \xi)^{\beta}
	e^{-2\pi i (k,\xi)} d\xi=\\
	(-1)^{[\beta]} \int\limits_{\rd} \h f(\xi) \overline{
	(2\pi i M^{*-j}\xi)^{\beta} e^{2\pi i (M^{-j}k,\xi)}} d\xi=
	(-1)^{[\beta]} m^{j/2}
	\langle
	\h f, \h {D^{\beta} \delta_{jk}}
	\rangle.
\ean
	If now $\w\phi = \sum\limits_{[\beta]\le N}
\overline{a_{\beta}} (-1)^{[\beta]} D^{\beta} \delta$ (we say that $\w\phi$
\textit{is associated with} $L$), then
	$$
m^{-j/2} L f(M^{-j}\cdot)(-k)=	
	\langle 	\h f, \h {\w\phi_{jk}} 	\rangle,\quad k\in\zd.
$$
Hence,
\be
Q_j(\phi,\w\phi, f)=m^{-j/2}\sum_{k\in\zd}L f(M^{-j}\cdot)(-k)\phi_{jk}.
\label{51}
\ee
We are interested  in the error of approximation 	of a function $f$ by these operators.
		
		\begin{theo}[\cite{KS1}]
\label{theoQjL}
	Let $2\le p \le \infty$, $1/p+1/q=1$,  a differential operator $L$ be defined by~(\ref{11}).
	Suppose
\begin{itemize}
  \setlength{\itemsep}{0cm}%
  \setlength{\parskip}{0cm}%

      \item[$(a)$]  $\w\phi$ is the distribution associated with $L$;

 	\item[$(b)$] $\phi \in {\cal L}_p $ and
	there exists $B_{\phi}>0$ such that
$
\sum_{k\in\zd}  |\h\phi(\xi+k)|^q<B_{\phi}$ for all $\xi\in\rd$;

	\item[$(c)$] there exist $n\in\n$ and $\delta\in(0, 1/2)$ such that
	$\h\phi\h{\w\phi}$ is  boundedly differentiable up to order $n$ on
	$\{|\xi|<\delta\}$,
	$\h\phi$ is boundedly	differentiable up to order $n$  on $\{|\xi+l|<\delta\}$	for all $l\in\zd\setminus\{\nul\}$;
	the function $\sum_{l\in\zd,\, l\neq\nul}|D^\beta \h \phi (\xi+l)|$
is bounded on	$\{|\xi|<\delta\}$  for $[\beta]=n$; $D^{\beta}(1-\h\phi\h{\w\phi})(0) = 0$ for $[\beta]<n$;
	the Strang-Fix condition of order $n$ holds for $\phi$;

\item[$(d)$]  $f\in L_p$,
       $\h f\in L_q$, $\h f(\xi)=O(|\xi|^{-N-d-\varepsilon})$
as $|\xi|\to\infty$,
	$\varepsilon>0$.
 	
\end{itemize}

\noindent
Then the following statements hold:
\begin{itemize}
  \setlength{\itemsep}{0cm}%
  \setlength{\parskip}{0cm}%

      \item[$(A)$]  if $n< N+\frac dp + \varepsilon$, then
$$
 \left\|f-m^{-j/2}\sum_{k\in\zd} L f(M^{-j}\cdot)(-k) \phi_{jk}\right\|_p\le C
 \vartheta^{-jn}
 $$
  for every  positive number $\vartheta$   which is smaller in module
than any eigenvalue of  $M$ and  some $C$ which does not depend on $j$;

 	\item[$(B)$] if $M$ is an isotropic matrix and $\lambda$
is its eigenvalue then
	 \be
	 \left\|f-m^{-j/2}\sum_{k\in\zd} L f(M^{-j}\cdot)(-k) \phi_{jk}\right\|_p\le
	 \begin{cases}
	 C |\lambda|^{-j(N+\frac dp + \varepsilon)}  &\mbox{if }
	n> N+\frac dp + \varepsilon\\
	  C (j+1)^{1/q} |\lambda|^{-jn} &\mbox{if }
	 n= N+\frac dp + \varepsilon \\
	C|\lambda|^{-jn}
	 &\mbox{if }
	 n< N+\frac dp + \varepsilon
	\end{cases},
	\label{fTheoQjL}
	\ee
	where $C$ does not depend on $j$;
	
	\item[$(C)$]  if  $\w\phi$ and $\phi$ are strictly compatible, then
	$$
 \left\|f-m^{-j/2}\sum_{k\in\zd} L f(M^{-j}\cdot)(-k) \phi_{jk}\right\|_p\le C
 \vartheta^{-j(N+\frac{d}{p}+\varepsilon)}
 $$
 for every  positive number $\vartheta$   which is smaller in module
than any eigenvalue of  $M$ and  some $C$ which does not depend on $j$.

\end{itemize}

\end{theo}	

Now we prove similar results for the class of  functions $\phi\in \cal B$.

\begin{theo}
\label{theoQjL_new}
		Let $2\le p < \infty$, $1/p+1/q=1$,
		a differential operator $L$ be defined by~(\ref{11}).
	Suppose
\begin{itemize}
  \setlength{\itemsep}{0cm}%
  \setlength{\parskip}{0cm}%

      \item[$(a)$]  $\w\phi$ is the distribution associated with $L$;

\item[$(b)$] $\phi\in \cal B$;

\item[$(c)$] $\w\phi$ and $\phi$ are strictly compatible;

\item[$(d)$]  $f\in L_p$,
       $\h f\in L_q$, $\h f(\xi)=O(|\xi|^{-N-d-\varepsilon})$
as $|\xi|\to\infty$,  $\varepsilon>0$.

\end{itemize}

\noindent
	Then
	\be
 \left\|f-m^{-j/2}\sum_{k\in\zd} L f(M^{-j}\cdot)(-k) \phi_{jk}\right\|_p\le C
 \vartheta^{-j(N+\frac{d}{p}+\varepsilon)},
 \label{12+}
 \ee
 where $C$ does not depend on $j$.
\end{theo}	
{\bf Proof.} 	Since, by~(\ref{51}),  all conditions of Theorem~\ref{theoQj_new}
are satisfied for $\w\phi$, $\phi$,
  (\ref{12+}) follows from~(\ref{8})
	and the next lemma.
$\Diamond$

\begin{lem}
Let $1\le q < \infty$, $1/p+1/q=1$, $N\in\z_+$, $\varepsilon>0$,
$g\in L_q$,
$g(\xi)=O(|\xi|^{-N-d-\varepsilon})$ as $|\xi|\to\infty$, $\gamma< N+\frac dp+\varepsilon$. Then
$$
	\|M^{*-j}\|^{\gamma q}  {\cal I}_{j,\gamma,q}^{\rm Out}(g) \le
	C \vartheta^{-j  q(N+\frac dp + \varepsilon)},
$$
	where $\vartheta$  is any positive number  which is smaller in module
than any eigenvalue of  $M$, $C$  does not depend on $j$.
	\label{lem2}
\end{lem}
{\bf Proof.} Throughout the proof we denote by $C_1, C_2,\dots$
different constants which do not depend on $j$.
Since
there exists  $A>0$ such that
$|g(\xi)|\le C_1 |\xi|^{-N-d-\varepsilon}$ for any $|\xi|>A$ and
  the set $\{|M^{*-j}\xi| \ge\delta\}$ is a subset of
	$\{|\xi|\ge\delta/\|M^{*-j}\|\}$, we have
		$$
	\|M^{*-j}\|^{\gamma q}  {\cal I}_{j,\gamma,q}^{\rm Out}(g)=
	 \|M^{*-j}\|^{\gamma q}
	\int\limits_{|M^{*-j}\xi|\ge\delta}
	|\xi|^{\gamma q}| g(\xi)|^q d\xi	\le
	C_1\|M^{*-j}\|^{\gamma q}
	\int\limits_{|\xi|\ge\delta/\|M^{*-j}\|}
	\frac {d\xi}
	{|\xi|^{(N+d+\varepsilon-\gamma)q}}
	$$
for all $j>j_0,$ where $j_0\in\z$ is such that 	$\frac {\delta} {\|M^{*-j_0}\|}>A.$
Using general polar coordinates with $\rho:=|\xi|$ and taking into account that
$(N+d+\varepsilon-\gamma)q>d$, we obtain
$$
\int\limits_{|\xi|\ge\delta/\|M^{*-j}\|}
	\frac {d\xi}
	{|\xi|^{(N+d+\varepsilon-\gamma)q}}\le
C_2\int\limits_{ {\delta}/{\|M^{*-j}\|}}^{+\infty}
	\frac {1}
	{\rho^{(N+d+\varepsilon-\gamma)q-d+1}} d\rho\le
C_3\|M^{*-j}\|^{q(N+\frac dp + \varepsilon-\gamma)},
$$
and, by~(\ref{00}),
$$
	\|M^{*-j}\|^{\gamma q}  {\cal I}_{j,\gamma,q}^{\rm Out}(g) \le	
C_4\|M^{*-j}\|^{q(N+\frac dp + \varepsilon)}\le
	C \vartheta^{-j  q(N+\frac dp + \varepsilon)}. \quad\Diamond
$$
\medskip

Theorem~\ref{theoQjL_new} says nothing about $p=\infty$. Now we consider this case and prove
a generalization of Brown's
inequality~(\ref{brown}).
\begin{theo}
		Let		a differential operator $L$ be defined by~(\ref{11}).
	Suppose
\begin{itemize}
  \setlength{\itemsep}{0cm}%
  \setlength{\parskip}{0cm}%

      \item[$(a)$]  $\w\phi$ is the distribution associated with $L$;

\item[$(b)$] $\phi\in \cal B$;

\item[$(c)$] $\w\phi$ and $\phi$ are strictly compatible;

\item[$(d)$]  $f\in C(\R^d)$,
       $\h f\in L$, $\h f(\xi)=O(|\xi|^{-N-d-\varepsilon})$
as $|\xi|\to\infty$,  $\varepsilon>0$.
\end{itemize}
Then for every $x\in\rd$ and $j\in \z$ the series
 $$
\sum_{k\in\zd} L f(M^{-j}\cdot)(-k) \phi_{jk}(x),
 $$
 considered as the limit of cubic partial sums, converges, and
\be
\left|f(x)- m^{-j/2}\sum_{k\in\zd} L f(M^{-j}\cdot)(-k) \phi_{jk}(x)\right|\le
C \|M^{*-j}\|^{N}\int\limits_{|M^{*-j}\xi|\ge \delta}|\xi|^N
 |\h f(\xi)|d\xi\le C' \vartheta^{-j(N+\varepsilon)},
\label{62}
\ee
where $C$ does not depend on $f$, $j$, and $x$; $C'$ does not depend on $j$ and $x$;
$\vartheta$  is a positive number  which is smaller in module
than any eigenvalue of  $M$.
\label{t8}
\end{theo}
{\bf Proof.}
Let $\phi$ be given by~(\ref{61}),
 $x\in\rd$. Set
$\Theta(\xi):=\sum_{s\in\,\zd}\theta(\xi+s) e^{2\pi i(x, \xi+s)}$.
By the Poisson summation formula, $\Theta$ is a summable $1$-periodic (with respect to
each variable) function and its $n$-th Fourier coefficient is
$$
\h\Theta(n)=
\int\limits_{\rd}
\theta(\xi)e^{-2\pi i(n-x,\xi)}\,d\xi=
\h\theta(n-x)=\phi(x-n).
$$
Since  $\Theta$ is a bounded function,  the cubic partial
Fourier sums are uniformly bounded in $L_\infty$-norm, and
the corresponding Fourier series
converges to the function almost everywhere.  Using this and
Lebesgue's dominated convergence theorem, for every $\beta\in\zd$, $[\beta]\le N$, we derive
\ban
\lim_{M\to\infty}\sum_{\|n\|_\infty\le M}D^\beta f(-n)\phi(x+n)=
\lim_{M\to\infty}
\int\limits_{\rd}\sum_{\|n\|_\infty\le M}\phi(x+n)e^{-2\pi i (n,\xi)}
(2\pi i\xi)^\beta\h f(\xi)\,d\xi=
\\
\int  \limits_{\rd}   \lim_{M\to\infty}  \sum_{\|n\|_\infty\le M}   \phi(x-n)e^{2\pi i (n,\xi)}
(2\pi i\xi)^\beta\h f(\xi)\,d\xi=
\int\limits_{\rd}
e^{2\pi i(x, \xi)} \sum_{s\in\,\zd}
\theta(\xi+s)e^{2\pi i(x, s)}(2\pi i\xi)^\beta\h f(\xi)\,d\xi.
\ean
Replacing $x$ by $M^jx$ and $f$ by
$f(M^{-j}\cdot)$, after a change of variable,  we obtain
\ban
\sum_{n\in\zd}D^\beta f(M^{-j}\cdot)(-n)\phi_{jn}(x)=
\lim_{M\to\infty}\sum_{\|n\|_\infty\le M}D^\beta f(M^{-j}\cdot)(-n)\phi_{jn}(x)=
\\
m^{j/2}\int\limits_{\rd}
e^{2\pi i(x, \xi)} \sum_{s\in\,\zd}\theta(M^{*-j}\xi+s)
e^{2\pi i (M^jx, s)}(2\pi iM^{*-j}\xi)^\beta \h f(\xi)\,d\xi.
\ean
This yields
\ba
m^{-j/2}\sum_{n\in\zd} Lf(M^{-j}\cdot)(-n)\phi_{jn}(x)=
\nonumber
\\
\int\limits_{\rd}
e^{2\pi i(x, \xi)} \sum_{s\in\,\zd}\theta(M^{*-j}\xi+s)
e^{2\pi i(M^jx, s)}\sum_{[\beta]\le N}\alpha_\beta(2\pi iM^{*-j}\xi)^\beta\h f(\xi)\,d\xi.
\label{45}
\ea
Set
$$
f_1(x)=\int\limits_{|M^{*-j}\xi|\le \delta}\h f(\xi)e^{2\pi i(x,\xi)}\,d\xi,\quad f_2(x)= f(x)-f_1(x).
$$
 If $|M^{*-j}\xi|\le \delta$, then
$$
\sum_{s\in\,\zd}\theta(M^{*-j}\xi+s)
e^{2\pi i (M^jx, s)}=\h\phi(M^{*-j}\xi)
$$
Taking into account that
$\h\phi(M^{*-j}\xi)\overline{\h{\w\phi}(M^{*-j}\xi)}=1$, we have
$$\sum_{s\in\,\zd}\theta(M^{*-j}\xi+s)
e^{2\pi i(M^jx, s)}\sum_{[\beta]\le N}\alpha_\beta(2\pi iM^{*-j}\xi)^\beta=1.
$$
Hence, it follows from~(\ref{45}) that
\be
m^{-j/2}\sum_{n\in\zd} Lf_1(M^{-j}\cdot)(-n)\phi_{jn}(x)=
\int\limits_{|M^{*-j}\xi|\le \delta}
e^{2\pi ix\xi} \h f(\xi)\,d\xi=f_1(x).
\label{46}
\ee
Since
$|e^{2\pi i(x,\xi)} \sum_{s\in\,\zd}\theta(M^{*-j}\xi+s)e^{2\pi i(M^jx,s)}|\le C_1$,
where $C_1$ depends only on $\phi$,
using~(\ref{45}) for $f_2$ and~(\ref{46}), we have
\ban
\left|f(x)- m^{-j/2}\sum_{k\in\zd} L f(M^{-j}\cdot)(-k) \phi_{jk}(x)\right|=
\left|f_2(x)- m^{-j/2}\sum_{k\in\zd} L f_2(M^{-j}\cdot)(-k) \phi_{jk}(x)\right|\le
\\
 \int\limits_{|M^{*-j}\xi|\ge\delta}
 \left(1+C_1\sum_{[\beta]\le N}|\alpha_l(2\pi M^{*-j}\xi)^\beta)|\right)|\h f(\xi)|d\xi\le
C \|M^{*-j}\|^{N}\int\limits_{|M^{*-j}\xi|\ge \delta}|\xi|^N
 |\h f(\xi)|d\xi,
\ean
which yields the first inequality in~(\ref{62}).
For  the second inequality it remains to apply Lemma~\ref{lem2} with $\gamma=N$, $q=1$. $\Diamond$

\begin{coro}
If, under  assumptions $(a)-(c)$ of Theorem~\ref{t8}, a function $f$ is such that
its Fourier transform is supported in $\{\xi:\ |M^{*-j}\xi|<\delta\}$,
where $\delta$ is from Definition~\ref{d1}, then
\be
f(x)=L f(M^{-j}\cdot)(-k) \phi_{jk}(x)\quad \text{for all}\quad x\in\rd.
\label{47}
\ee
If, moreover, $\h\phi$ is supported in $[-1/2,1/2]^d$ and
$\h\phi(\xi)\overline{\h{\w\phi}(\xi)}=1$ for all $\xi\in [-1/2,1/2]^d$, then~(\ref{47})
holds for every $f$ whose Fourier transform is supported in $M^{*j}[-1/2, 1/2]^d$.
\label{coro1}
\end{coro}	

{\bf Proof.}
The first statement is a trivial consequence of Theorem~\ref{t8}. Analysing the proof of this theorem, it is easy to see that the second statement is also true.
$\Diamond$

\bigskip

Corollary~\ref{coro1} is an analog of the classical sampling theorem.	

Now consider the following differential equation
$$
Lf:=\sum_{[\beta]\le N} a_\beta D^{\beta} f =g,
$$
where $a_{\nul}\ne0$ and $g$ is a function band-limited to $[-1/2, 1/2]^d$.
Let $\w\phi$ is a distribution associated with $L$. Assume that $\h{\w\phi}$
does not vanish on $[-1/2, 1/2]^d$ and define $\phi$ by
$$
\h\phi(\xi)=
 \begin{cases}
	 	\frac1{\h{\w\phi}(x)}, &\mbox{if } \xi\in [-1/2, 1/2]^d,
\\
	0, &\mbox{if } \xi\not\in [-1/2, 1/2]^d.
	\end{cases}
$$
By Corollary~\ref{coro1}, we have
$$
f(x)=\sum_{k\in\zd} g(k)\phi(x+k).
$$
For example, if $d=1$, $Lf=f-f''$, then $\phi(x)=\int_{-\pi}^{\pi}\frac{\cos tx\,dt}{1+t^2}$.
The latter function can be expressed via four special functions: hyperbolic sine, hyperbolic cosine,
 integral sine and  integral cosine, see Wolfram  Mathematica.

\section{Falsified sampling expansions}
\label{FalSamExp}

 If exact sampled values of a signal $f$ are  known,
 then sampling expansions are very useful for applications.
 Theorems~\ref{theoQjL},  \ref{theoQjL_new}, and \ref{t8}  provide error estimates
 for this case.
  Now we discuss what happens
 if exact sampled values  $f(-M^{-j}k)$
  are replaced by  average values. We assume that at each point $M^{-j}k$
	one knows the following average value of $f$
\be
	\frac  1 {V_{h(u)}} \int\limits_{B_{h(u)}} f(M^{-j}k+M^{-j}t)d t=
	\frac  {m^j} {V_{h(u)}} \int\limits_{M^{-j}B_{h(u)}} f(M^{-j}k+t)d t=:
{\rm Av}_{h{(u)}}(f, M^{-j}k),
	\label{fFAverValue}
	\ee
where $h(u)$ is a positive function defined on $(0, \infty)$
and $u$ is a random value with probability density~$w$;
$V_h$ is the volume of the ball $B_h$.
Set
$$
E(f, M^{-j}k)=E(f, M^{-j}k, h,w)=\int\limits_0^\infty du\,w(u){\rm  Av}_{h(u)}(f, M^{-j}k).
$$
Note that if $h(u)\equiv h>0$, then $E(f, M^{-j}k)={\rm Av}_{h}(f, M^{-j}k)$.

We  are interested in error analysis for the {\em falsified sampling expansions}
$$
m^{-j/2}\sum_{k\in\zd}E(f, M^{-j}k)\,\phi_{jk}.
$$
First we prove the following lemma.
\begin{lem}
Let $N\in\n$, let a function $f$
be continuously differentiable up to order $N$, and let $A$ be a real $d\times d$ matrix.
	Then for all $t,x\in\rd$
	$$
	\sum\limits_{[\beta]<{N+1}}
	\frac{D^{\beta} f (A x)} {\beta!}
	(A t)^{\beta}=
	\sum\limits_{[\beta]<{N+1}}
	\frac{ D^{\beta} f(A\cdot)(x)}{\beta!}
	t^{\beta}.
	$$
\label{lemDiffMatrixA}
\end{lem}

{\bf Proof. }
First we introduce some additional notations. Let $r\in\z_+,$
$O_r=\{\beta\in\zd_+:[\beta]=r\}$. Assume that the set $O_r$ is
ordered by lexicographic order. Namely,
$(\beta_1,\dots,\beta_d)$ is less than
$(\alpha_1,\dots,\alpha_d)$ in lexicographic order
if $\beta_j=\alpha_j$ for $j=1,\dots,i-1$ and
$\beta_i<\alpha_i$ for some $i$.
Let $S(A,r)$ be a  $(\# O_r)\times (\# O_r)$ matrix
 which is uniquely determined by
	\be
	\frac {(A t)^{\alpha}}{\alpha!} = \sum_{\beta\in O_r}
	[S(A,r)]_{\alpha,\beta} \frac {t^{\beta}}{\beta!},
	\label{fBHanEq}	
	\ee
where $\alpha\in O_r,$ $t\in\rd.$ It can be verified that
	\be
	\alpha! [S(A,r)]_{\alpha,\beta}=
	\beta! [S(A^*,r)]_{\beta,\alpha}.
	\label{fBHanEq2}\ee
The above notation and the latter fact is borrowed from~\cite{HanSymSmooth}.

Fix $\beta \in\zd_+.$ Let ${\cal E}$ be the set of ordered samples with replacement
of size $[\beta]$ from the set $\{e_1,\dots,e_d\}$, where
$e_k$ is the $k$-th ort in $\rd$. An element $e\in {\cal E} $ is a set  $\{e_{i_1},\dots,e_{i_{[\beta]}}\}$,
where $i_l\in\{1,\dots,d\},$ $l=1,\dots,[\beta],$
$\#{\cal E} = d^{[\beta]}.$ For $e\in {\cal E}$ denote by $(e)_l:=e_{i_l}$,  $l=1,\dots,[\beta].$
Let $T$ be a function defined on ${\cal E}$ by $T(e):=\sum_{i=1}^{[\beta]} (e)_i.$
Note that $T(e)\in \zd_+$ and $[T(e)]=[\beta].$ Denote by $b$ an element of
${\cal E}$ so that $T(b)=\beta$. Such $b$ is unique up to a permutation.
Using the higher chain rule,  we have
	$$
	D^{\beta} f(A\cdot)(x)= \frac {\partial ^{[\beta]} f(A\cdot)}
	{\partial x^{\beta}}(x)=
	\sum_{e\in{\cal E}} \frac {\partial^{[\beta]} f(y)}
	{\partial y^{T(e)}}\Big|_{y=Ax}
	\prod_{i=1}^{[\beta]} \frac {\partial (Ax)^{(e)_i}}
{\partial x^{(b)_i}},\quad x\in\rd,
	$$
	where $\prod_{i=1}^{[\beta]} \frac {\partial (Ax)^{(e)_i}}{\partial x^{(b)_i}}$	 does not depend on $x$ and $D^{\beta} f(A\cdot)(x)$ does not depend on the choice of $b$.
For different elements $e,h\in{\cal E}$, we may have $T(e)=T(h).$
Thus, we can group terms in the sum with equal values of $T(\cdot)$. Namely,
\be
	D^{\beta} f(A\cdot)(x) =
		\sum_{\alpha \in\zd_+, [\alpha]=[\beta]}
		\frac {\partial^{[\beta]} f(y)}
	{\partial y^{\alpha}}\Big|_{y=Ax}
	\sum_{e\in{\cal E}, T(e)=\alpha}
	\prod_{i=1}^{[\beta]} \frac {\partial (Ax)^{(e)_i}}
	{\partial x^{(b)_i}}.
\label{16}
\ee	
If $f(x)=e^{2\pi i (t,x)},$  $t\in\rd,$ then
	$
D^{\beta} f(A\cdot)(x) =	D^{\beta} e^{2\pi i (A^*t,x)}=(A^*t)^{\beta} e^{2\pi i (t,Ax)}.
	$
On the other hand, by~(\ref{16}),
	$$
	D^{\beta} f(A\cdot)(x) =
			 e^{2\pi i (t,Ax)}
			 \sum_{\alpha \in\zd_+, [\alpha]=[\beta]}
t^{\alpha}
	\sum_{e\in{\cal E}, T(e)=\alpha}
	\prod\limits_{i=1}^{[\beta]}
	\frac {\partial (Ax)^{(e)_i}}{\partial x^{(b)_i}}.
	$$
	Thus, due to~(\ref{fBHanEq}) with the matrix $A$ replaced by $A^*$,
	 we obtain
	\be\sum_{e\in{\cal E}, T(e)=\alpha}
	\prod\limits_{i=1}^{[\beta]}
	\frac {\partial (Ax)^{(e)_i}}{\partial x^{(b)_i}}=[S(A^*,r)]_{\beta,\alpha} \frac {\beta!}{\alpha!}.
	\label{fLemBHDiff}	
	\ee
	
Now let $f$ be an arbitrary  function continuously differentiable  up
	to order $N,$ $0\le r \le N$, $r\in\z_+.$
	It follows from~(\ref{fBHanEq2}) 	and~(\ref{fLemBHDiff}) that
	
	$$
	\sum_{\alpha\in O_r}
	\frac {D^{\alpha} f(A x)}{\alpha!} (A t)^{\alpha}=
		\sum_{\alpha\in O_r}
	D^{\alpha} f(A x) \sum_{\beta\in O_r}
[S(A,r)]_{\alpha,\beta} \frac {t^{\beta}}{\beta!}=
		\sum_{\beta\in O_r}  t^{\beta}
		\sum_{\alpha\in O_r}
	D^{\alpha} f(A x)
 \frac {[S(A^*,r)]_{\beta,\alpha}}{\alpha!}=
	$$
	$$
	\sum_{\beta \in O_r}  \frac {t^{\beta}}{\beta!}
		\sum_{\alpha\in O_r}
	D^{\alpha} f(A x)
 \sum_{e\in{\cal E}, T(e)=\alpha}
	\prod\limits_{i=1}^{[\beta]} \frac {\partial (Ax)^{(e)_i}}
	{\partial x^{(b)_i}}=
	\sum_{\beta\in O_r}  \frac {t^{\beta}}{\beta!}
	D^{\beta} f(A\cdot)(x).
	$$
It remains to sum the latter expression over $r$ from $0$ to $N.$
$\Diamond$
\medskip

Let
	\be
	L f(M^{-j}\cdot)(k)=
	\sum_{[\beta]<{N+1}} a_\beta D^{\beta}
	f(M^{-j}\cdot)(k) ,
	\quad
	a_\beta =\int\limits_0^\infty du\frac  {w(u)} {\beta! V_{h(u)}} \int\limits_{B_{h(u)}} t^{\beta}d t.
	\label{fOperL}
	\ee
By Lemma~\ref{lemDiffMatrixA}, we have
\begin{equation}\label{15}
  \begin{split}
    	L f(M^{-j}\cdot)(k)=
	 \int\limits_0^\infty du\frac  {w(u)} { V_{h(u)}} \int\limits_{B_{h(u)}}
	 \sum_{[\beta]<{N+1}}	
	\frac {D^{\beta}
	f(M^{-j}\cdot)(k)}{\beta!}
	t^{\beta}d t &=
	 	\\
\int\limits_0^\infty du\frac  {w(u)} { V_{h(u)}} \int\limits_{B_{h(u)}}
	 \sum_{[\beta]<{N+1}}	
	\frac{D^{\beta} f (M^{-j} k)} {\beta!}
	(M^{-j} t)^{\beta} d t &=
\\
{m^j}\int\limits_0^\infty du\frac  {w(u)}{V_{h(u)}}
	\int\limits_{M^{-j}B_{h(u)}}
	 \sum_{[\beta]<{N+1}}	
	\frac{D^{\beta} f (M^{-j} k)} {\beta!} 		t^{\beta}d t&,
  \end{split}
\end{equation}
and, due to the Tailor formula,
$$
  {m^{j}}\int\limits_0^\infty du\frac  {w(u)} {V_{h(u)}}
	\int\limits_{M^{-j}B_{h(u)}} f(M^{-j}k+t)d t
	\approx L f(M^{-j})(k).
$$
To investigate the convergence and  approximation order
of falsified sampling expansions we can use
Theorems~\ref{theoQjL}, \ref{theoQjL_new}, \ref{t8},   and estimate the sum
$
m^{-j/2}\sum_{k\in\zd}\varepsilon_j(-k)\,\phi_{jk},
$
where
	\ba
	\varepsilon_j(k):=\int\limits_0^\infty w(u){\rm  Av}_{h{(u)}}(f, M^{-j}k)\,du-
	L f(M^{-j}\cdot)(k).
		\label{fEpsK}
	\ea

\begin{prop}
Let $d< p \le \infty$,   $\phi \in {\cal L}_p$ or  $\phi\in \cal B$, $p\neq\infty$.
Suppose  $N\in\n$, $f\in W_p^{N+1}$, the operator $L$ is defined by~(\ref{fOperL}),
$\varepsilon_j(k)$ is defined by~(\ref{fEpsK}), $w$ and $h$ are as in~(\ref{fFAverValue}). Then
\be
	\left\| m^{-j/2}
	\sum\limits_{k\in\zd}
	\varepsilon_j(-k) \phi_{jk}\right\|_{p} \le
	C \int\limits_0^\infty w(u)(1+h^{N+1}(u))\,du\,\|f\|_{W_p^{N+1}}\,\vartheta^{-j(N+1)}
\label{13}
\ee
for every  positive number $\vartheta$   which is smaller in module
than any eigenvalue of  $M$ and  some  $C$ which  does not depend
on $f$, $h$, $w$, and $j$.
\label{t1}
\end{prop}

{\bf Proof.} Let us fix $j\in\n,$ $\varepsilon_j=\{\varepsilon_j(k)\}_{k\in\zd}$.	
If $\phi \in {\cal L}_p$, then, due to Proposition~\ref{propLp},
	\be
	\left\| m^{-j/2}	\sum\limits_{k\in\zd}
	\varepsilon_j(-k) \phi_{jk}\right\|_{p} =
	m^{-j/p} \left\| \sum\limits_{k\in\zd}
	\varepsilon_j(-k)\phi_{0k}\right\|_p\le
m^{-j/p}\|\phi\|_{{\cal L}_p}\|\varepsilon_j\|_{\ell_p}.
	\label{14}
\ee
If $\phi \in {\cal B}$ and $p\ne\infty$, then, due to Proposition~\ref{prop2},
	\be
	\left\| m^{-j/2}	\sum\limits_{k\in\zd}
	\varepsilon_j(-k) \phi_{jk}\right\|_{p} =
	m^{-j/p} \left\| \sum\limits_{k\in\zd}
	\varepsilon_j(-k)\phi_{0k}\right\|_p\le
m^{-j/p} C_{\phi, q}\|\varepsilon_j\|_{\ell_p}.
	\label{14'}
\ee

By the Taylor formula with 	integral remainder, we have
	$$
	f(M^{-j}k+t)= \sum\limits_{[\beta]<{N+1}}
	\frac {D^{\beta} f(M^{-j}k)}{\beta!} t^{\beta} +
	\sum\limits_{\beta\in\zd_+, [\beta]=N+1}
	\frac {N+1}{\beta!} t^{\beta}
	\int\limits_0^1 (1-\tau)^N D^{\beta} f(M^{-j}k+t\tau)d\tau.
	$$
It follows from~(\ref{15})  that
$$
\varepsilon_j(k)={m^j} \int\limits_0^\infty du\frac {w(u)} {V_{h(u)}} \int\limits_{M^{-j}B_{h(u)}}\,dt
\sum\limits_{\beta\in\zd_+, [\beta]=N+1}
	\frac {N+1}{\beta!} t^{\beta}
	\int\limits_0^1 (1-\tau)^N D^{\beta} f(M^{-j}k+t\tau)\,d\tau.
$$
	Hence, taking into account that $|t^{\beta}|\le |t|^{[\beta]}$, we get
$$
	|\varepsilon_j(k)|\le
	\sum\limits_{\beta\in\zd_+, [\beta]=N+1}
	\frac {N+1}{\beta!}
	{m^j} \int\limits_0^\infty du\frac {w(u)} {V_{h(u)}}\int\limits_{M^{-j}B_{h(u)}} dt\,|t|^{[\beta]}
	\int\limits_0^1  |D^{\beta} f(M^{-j}k+t\tau)|d\tau \le
	$$
		$$
	\sum\limits_{\beta\in\zd_+, [\beta]=N+1}
	\frac {N+1}{\beta!} {m^j} \int\limits_0^\infty du\frac {w(u)} {V_{h(u)}}	
	\int\limits_0^1   d\tau	 \int\limits_{\tau M^{-j}B_{h(u)}}
	\frac {|t|^{[\beta]} }{\tau^{[\beta]}}
	| D^{\beta} f(M^{-j}k+t) |\frac {dt}{\tau^d} =
	 $$
	 $$
	 \sum\limits_{\beta\in\zd_+, [\beta]=N+1}
	 \frac {N+1}{\beta!} \int\limits_0^\infty du
	\frac {w(u)} {V_{h(u)}}	\int\limits_0^1   d\tau
	 \int\limits_{\tau B_{h(u)}}|M^{-j}t|^{[\beta]}
	| D^{\beta} f(M^{-j}k+M^{-j}t)| \frac {dt}{\tau^{[\beta]+d}}.
	 $$
 The latter integration is taken over the set
 $\{\tau\in[0,1], t\in \tau B_{h(u)}\}=\{\tau\in[0,1], |t|\le \tau {h(u)} \}$, or equivalently
$\{|t|\in[0,{h(u)}], \frac {|t|}{h(u)}\le\tau \le 1\}$. Changing the order of integration,
we obtain
	 $$|\varepsilon_j(k)|\le
	 \sum\limits_{\beta\in\zd_+, [\beta]=N+1}
	 \frac {N+1}{\beta!}\int\limits_0^\infty du\frac {w(u)} {V_{h(u)}}
	\int\limits_{B_{h(u)}} |M^{-j}t|^{[\beta]}
	 |D^{\beta} f(M^{-j}k+M^{-j}t)| dt \int\limits_{\frac {|t|}{h(u)}}^1
	\frac {d\tau}  {\tau^{N+d+1}} \le
	$$
	$$
	\sum\limits_{\beta\in\zd_+, [\beta]=N+1}
	 \frac {(N+1)}{(N+d)\beta!} \int\limits_0^\infty du\frac {w(u)} {V_{h(u)}}
	\int\limits_{B_{h(u)}}
	 |M^{-j}t|^{N+1}
	| D^{\beta} f(M^{-j}k+M^{-j}t) |
	  \left(\frac {h(u)}{|t|}\right)^{N+d} dt \le
	  $$
	  $$
	 	\sum\limits_{\beta\in\zd_+, [\beta]=N+1}
	 \frac {\|M^{-j}\|^{N+1}}{\beta! }
	 \int\limits_0^\infty du\frac {w(u)h^{N+d}(u)} {V_{h(u)}}
	\int\limits_{B_{h(u)}}		\frac{| D^{\beta} f(M^{-j}k+M^{-j}t) |}{|t|^{d-1}}
	   dt =
$$
$$
\sum\limits_{\beta\in\zd_+, [\beta]=N+1}
	 \frac {\|M^{-j}\|^{N+1} {m^j}}{\beta! }
	\int\limits_0^\infty du\frac {w(u)h^{N}(u)} {V_1}
	\int\limits_{M^{-j}B_{h(u)}}\frac{| D^{\beta} f(M^{-j}k+t) |}{|M^j t|^{d-1}}dt.
	  $$
If $p=\infty$, then
$$
\int\limits_{M^{-j}B_{h(u)}}\frac{| D^{\beta} f(M^{-j}k+t) |}{|M^j t|^{d-1}}dt\le
\int\limits_{M^{-j}B_{h(u)}}\frac{dt}{|M^j t|^{d-1}}\|f\|_{W_\infty^{N+1}}.
$$
If  $p\neq\infty$, then using H\"older's inequality, we have	
$$
\int\limits_{M^{-j}B_{h(u)}}\frac{| D^{\beta} f(M^{-j}k+t) |}{|M^j t|^{d-1}}dt\le
\lll\int\limits_{M^{-j}B_{h(u)}}{| D^{\beta} f(M^{-j}k+t) |^p}dt\rrr^{1/p}
\lll\int\limits_{M^{-j}B_{h(u)}}\frac{dt}{|M^j t|^{q(d-1)}}\rrr^{1/q},
$$
where $q=\frac p{p-1}$.
Since
$$
\int\limits_{M^{-j}B_{h(u)}}\frac{dt}{|M^j t|^{q(d-1)}}= m^{-j}
\int\limits_{B_{h(u)}}\frac{dt}{|t|^{q(d-1)}}=
m^{-j}h^{-q(d-1)+d}(u)
\int\limits_{B_1}\frac{dt}{|t|^{q(d-1)}}
$$
and $q(d-1)<d$, the latter integral is finite.
Summarizing the above estimates, we obtain
$$
|\varepsilon_j(k)|^p\le C_1 m^{j}
\|M^{-j}\|^{p(N+1)}\int\limits_0^\infty du\, w(u)h^{p(N-d+1+d/q)}(u)
\int\limits_{M^{-j}k+M^{-j}B_{h(u)}}{\sum\limits_{\beta\in\zd_+, [\beta]=N+1}| D^{\beta} f(t) |^p}dt
$$
and
$$
\sum_{k\in\zd}|\varepsilon_j(k)|^p\le C_2 m^{j}
\|M^{-j}\|^{p(N+1)}\int\limits_0^\infty du\, w(u)h^{p(N-d+1+d/q)}(u)
(1+h^{d}(u))\|f\|_{W_p^{N+1}}^p,
$$
where $C_2$ does not depend on $f$, $h$, and $j.$
It follows that
$$
\sum_{k\in\zd}|\varepsilon_j(k)|^p\le C_2  m^{j} \|M^{-j}\|^{p(N+1)}
\|f\|_{W_p^{N+1}}^p\int\limits_0^\infty  w(u)(h^{p(N-d+1+d/q)}(u)+h^{p(N+1)}(u))\, du.
$$ 
Combining this with~(\ref{14}), ~(\ref{14'}) and~(\ref{00}), we get~(\ref{13}).  $\Diamond$

\begin{prop}
Let $d< p' <\infty$,    $\phi\in \cal B$.
Suppose  $N\in\n$, $f\in W_{p'}^{N+1}$, operator $L$ is defined by~(\ref{fOperL}),
$\varepsilon_j(k)$ is defined by~(\ref{fEpsK}), $w$ and $h$ are as in~(\ref{fFAverValue}). Then
\be
	\left\| m^{-j/2}
	\sum\limits_{k\in\zd}
	\varepsilon_j(-k) \phi_{jk}\right\|_{\infty} \le
	C \int\limits_0^\infty w(u)(1+h^{N+1}(u))\,du\,\|f\|_{W_{p'}^{N+1}}\,\vartheta^{-j(N+1)}
\label{13'}
\ee
for every  positive number $\vartheta$   which is smaller in module
than any eigenvalue of  $M$ and  some  $C$ which  does not depend
on $f$, $h$, $w$ and $j$.
\label{t1'}
\end{prop}

{\bf Proof.} Let us fix $j\in\n,$ $\varepsilon_j=\{\varepsilon_j(k)\}_{k\in\zd}$.	
Because of~(\ref{eqK7.5}), it is easy to see that
$$
\sum_{k\in \zd}|\phi(x+k)|^q< C'_{\phi, q}
$$
for every $q>1$ and every $x\in\rd$. It follows that
	$$
	\left\| m^{-j/2}	\sum\limits_{k\in\zd}
	\varepsilon_j(-k) \phi_{jk}\right\|_{\infty}
	=\sup\limits_{x\in\rd} \left|\sum\limits_{k\in\zd}\varepsilon_j(-k) \phi(M^jx+k)\right|\le
	 (C^\prime_{\phi, q'})^{1/q'}\|\varepsilon_j\|_{\ell_{p'}},
	$$
where $1/q'+1/p'=1$. To get~(\ref{13'}) it remains to repeat all arguments of the proof of
Theorem~\ref{t1} after relation~(\ref{14'}), replacing $p$ and $q$ by $p'$ and $q'$ respectively.$\Diamond$

\begin{rem}
If  $M$ is an isotropic matrix
for which $\lambda$  is an eigenvalue,
then in the proof of Propositions~\ref{t1}, \ref{t1'} we can use inequality~(\ref{10}) instead of~(\ref{00}).
Hence $\vartheta$ in~(\ref{13}) and~(\ref{13'})   can be replaced by~$\lambda$.
\label{rem+}
\end{rem}

Using the above results we can state the convergence and  approximation order of
falsified sampling expansions.

\begin{theo}
\label{theoQjAverage}
Let   $d< p \le \infty$, $p\ge2$,  $1/p+1/q=1$,  $N\in\z_+,$
$M$ be an isotropic matrix dilation and $\lambda$ be its eigenvalue.
Suppose $\phi$ and $n$ are as
of Theorem~\ref{theoQjL}, where $\w\phi$ is the distribution associated with the differential
operator $L$ given by~(\ref{fOperL}), $w$ and $h$ are as in~(\ref{fFAverValue}) and
$$
\int\limits_0^\infty w(u)(1+h^{N+1}(u))\,du<\infty;
$$
$f\in L_p$,
       $\h f\in L_q$, $\h f(\xi)=O(|\xi|^{-N-1-\frac dq-\varepsilon})$
as $|\xi|\to\infty$, $\varepsilon>0$.
Then
\be
\Big\|f-m^{-j/2}\sum_{k\in\zd} E(f, M^{-j}k)\,\phi_{jk}\Big\|_p\le
 \begin{cases}
	 	C |\lambda|^{-j(N+1)}  &\mbox{if }
	
	n> N+1
\\
	C|\lambda|^{-jn}
	 &\mbox{if }
	 n\le N+1
	 	 \\
	\end{cases},
	 \label{fTheoQj3}
	\ee
where $C$ does not depend on  $j$.
\end{theo}

\textbf{{Proof.}} The proof follows from Theorem~\ref{theoQjL}, item (B), Proposition~\ref{t1}, and Remark~\ref{rem+}. $\Diamond$

\begin{theo}
\label{theoQjAverage2}
Let   $d< p \le \infty$, $p\ge2$, $1/p+1/q=1$, $N\in\z_+.$
Suppose  $\phi$ is as in item $(b)$ of Theorem~\ref{theoQjL} or $\phi\in{\cal B}$, $p\neq \infty$;
$\w\phi$ is the distribution associated with the differential
operator $L$ given by~(\ref{fOperL}),
$\w\phi$ and $\phi$ are strictly compatible;
$w$ and $h$ are as in Theorem~\ref{theoQjAverage}; $f\in L_p$,   $\h f\in L_q$,
 $\h f(\xi)=O(|\xi|^{-N-1-\frac dq-\varepsilon})$ as $|\xi|\to\infty$, $\varepsilon>0$.
Then
\be
\Big\|f-m^{-j/2}\sum_{k\in\zd}E(f, M^{-j}k)\,\phi_{jk}\Big\|_p\le C \vartheta^{-j(N+1)}
 	 \label{fTheoQj32_new}
	\ee
for every  positive number $\vartheta$   which is smaller in module
than any eigenvalue of  $M$ and  some  $C$ which  does not depend on $j$.
\end{theo}

\textbf{{Proof.}} The proof follows from Theorem~\ref{theoQjL}, item (C),
Theorem~\ref{theoQjL_new}, and Proposition~\ref{t1}. $\Diamond$

\bigskip

\begin{rem}
If the assumption on $f$ in Theorems~\ref{theoQjAverage} and~\ref{theoQjAverage2} is replaced by $f\in W_p^{N+1}$, $\h f\in L_q$,
 $\h f(\xi)=O(|\xi|^{-N-1-\frac dq})$ as $|\xi|\to\infty$, then  Theorem~\ref{theoQjAverage2} remains to be true, and inequality~(\ref{fTheoQj3}) in Theorem~\ref{theoQjAverage} must be replaced by
 $$
\Big\|f-m^{-j/2}\sum_{k\in\zd} E(f, M^{-j}k)\,\phi_{jk}\Big\|_p\le
 \begin{cases}
	 	C |\lambda|^{-j(N+1)}  &\mbox{if }
	
	n> N+1
\\
C (j+1)^\frac1q|\lambda|^{-j(N+1)}  &\mbox{if }  n=N+1
\\
	C|\lambda|^{-jn}
	 &\mbox{if }
	 n\le N+1
	 	 \\
	\end{cases}.
	$$
\end{rem}

\begin{theo}
\label{theoQjAverage2'}
Let   $d< p' < \infty$,  $1/p'+1/q'=1$, $N\in\z_+.$
Suppose  $\phi\in{\cal B}$,
$\w\phi$ is the distribution associated with the differential
operator $L$ given by~(\ref{fOperL}),
$\w\phi$ and $\phi$ are strictly compatible;
$w$ and $h$ are as in Theorem~\ref{theoQjAverage}; $f\in L_{p'}$,   $\h f\in L_{q'}$,
 $\h f(\xi)=O(|\xi|^{-N-1-\frac d{q'}-\varepsilon})$ as $|\xi|\to\infty$, $\varepsilon>0$.
Then
\be
\Big\|f-m^{-j/2}\sum_{k\in\zd}E(f, M^{-j}k)\,\phi_{jk}\Big\|_\infty\le C \vartheta^{-j(N+1)}
 	 \label{fTheoQj32_new'}
	\ee
for every  positive number $\vartheta$   which is smaller in module
than any eigenvalue of  $M$ and  some  $C$ which  does not depend on $j$.
\end{theo}

\textbf{{Proof.}} The proof follows from Theorem~\ref{t8}  and Proposition~\ref{t1'}. $\Diamond$

\bigskip


\begin{rem}
If the assumption on $f$ in Theorem~\ref{theoQjAverage2'} is replaced by $f\in W_{p'}^{N+1}$, $\h f\in L_{1}$,
$\h f(\xi)=O(|\xi|^{-N-d-\varepsilon})$ as $|\xi|\to\infty$, $\varepsilon>0$, then  Theorem~\ref{theoQjAverage2} remains to be true with
$$
\Big\|f-m^{-j/2}\sum_{k\in\zd}E(f, M^{-j}k)\,\phi_{jk}\Big\|_\infty\le C \vartheta^{-j(N+\min(1,\varepsilon))}
$$
instead of~\ref{fTheoQj32_new'}.
\end{rem}

Observing the proof of Proposition~\ref{t1}, one can see that in the
one-dimensional case  an analog of~(\ref{13})
holds  true for    a wider class of functions $f$ and any $p\ge1$.
Indeed, in this case  we have
$$
|\varepsilon_j(k)|\le
	 \frac {1}{(N+1)! }\int\limits_0^\infty du\frac {w(u)} {V_{h(u)}}
	\int\limits_{B_{h(u)}}
	 |M^{-j}t|^{N+1}
	|  f^{(N+1)}(M^{-j}k+M^{-j}t) |
	  \left(\frac {h(u)}{|t|}\right)^{N+1} dt \le
$$
$$
\frac {M^{-j(N+1)}}{(N+1)! }
	\int\limits_0^\infty du\frac {w(u)h^{N+1}(u)} {V_{h(u)}}\int\limits_{B_{h(u)}}	
		|  f^{(N+1)}(M^{-j}k+M^{-j}t) |	   dt =
	$$
	$$
\frac {M^{-jN}}{(N+1)! }\int\limits_0^\infty du\frac {w(u)h^{N+1}(u)} {V_{h(u)}}
	\int\limits_{M^{-j}B_{h(u)}+M^{-j}k}	|  f^{(N+1)}(t) |dt.
$$
It follows that
$$
\left\| m^{-j/2}	\sum\limits_{k\in\z}
	\varepsilon_j(-k) \phi_{jk}\right\|_{p} \le
	m^{-j/p}\|\phi\|_p
\sum_{k\in\z}|\varepsilon_j(k)|\le C\int\limits_0^\infty  w(u)h^{N}(u)du \|f\|_{W_{1}^{N+1}}M^{-j(N+\frac1p)}.
$$
This yields the following statements.

\begin{prop}
Let $d=1$, $p\ge1$, $N\in\n$, $\phi\in L_p$.
Suppose  $f\in W_1^{N+1}$,
$\varepsilon_j(k)$ is defined by~(\ref{fEpsK}), $w$ and $h$ are as in~(\ref{fFAverValue}) and
$$
\int\limits_0^\infty w(u)(1+h^{N}(u))\,du<\infty.
$$ Then
\be
	\left\| M^{-j/2}
	\sum\limits_{k\in\z}
	\varepsilon_j(-k) \phi_{jk}\right\|_{p} \le  C  \|f\|_{W_{1}^{N+1}} M^{-j(N+\frac 1p)},
\label{17}
\ee
where  $C$ does not depend on $f$ and $j$.
\label{t11}
\end{prop}

\begin{theo}
\label{theoQjAverage1}
Let $d=1$,  $2\le p \le \infty$,    $N\in\z_+$.
Suppose $\phi$ and $n$ are as in
of Theorem~\ref{theoQjL}, where $\w\phi$ is the distribution associated with the differential
operator $L$ given by~(\ref{fOperL}), or $\phi\in{\cal B}$;
$w$ and $h$ are as in Proposition~\ref{t11};
 $f\in W_1^{N+1}$, $f^{(N+1)}\in\ {\rm Lip}_{L_1} \,\varepsilon$,
$\varepsilon>0$. Then
\be
\Big\|f-M^{-j/2}\sum_{k\in\z} E(f, M^{-j}k)\,\phi_{jk}\Big\|_p\le
 \begin{cases}
	 	C M^{-j(N+\frac 1p )}  &\mbox{if }
	
	n> N+\frac 1p
\\
	C M^{-jn}
	 &\mbox{if }
	 n\le N+\frac 1p
	 	 \\
	\end{cases},
	 \label{fTheoQj31}
	\ee
where $C$ does not depend on $j$.
\end{theo}

\begin{theo}
Let $d=1$,  $2\le p \le \infty$,   $N\in\z_+$.
Suppose  $\phi$ is as in item $(b)$ of Theorem~\ref{theoQjL} or $\phi\in{\cal B}$;
$\w\phi$ is the distribution associated with the differential
operator $L$ given by~(\ref{fOperL}),
$\w\phi$ and $\phi$ are strictly compatible; $w$ and $h$ are as in Proposition~\ref{t11};
$f\in W_1^{N+1}$, $f^{(N+1)}\in\ {\rm Lip}_{L_1} \,\varepsilon$,
$\varepsilon>0$.
Then
\be
\Big\|f-M^{-j/2}\sum_{k\in\z} E (f, M^{-j}k)\,\phi_{jk}\Big\|_p\le C M^{-j(N+\frac1p)}
 	 \label{fTheoQj32}
	\ee
for  some  $C$ which  does not depend on $j$.
\end{theo}

	\section{Examples}
	\label{example}

In this section some examples will be given to illustrate the above results.

I. First we discuss construction of band-limited functions $\phi$.
For every differential operator $L$ one can easily
construct  $\phi$  supported on a small neighborhood of zero and such that
on the same neighborhood $\h\phi \h{\w\phi}=1$, where $\w\phi$ is a distribution associated with
$L$.   The function $\phi$ is in $\cal B$,
$\w\phi$ and $\phi$ are strictly compatible,
so, due to Theorems~\ref{theoQjL_new} and \ref{theoQjAverage2},
the  approximation order of the corresponding differential and falsified
expansions with arbitrary matrix dilation depends only on  how smooth is a function $f$.
In particular, in the case $Lf=f$,  the function
$\phi(x)=\prod_{k=1}^d \frac {\sin  {\pi x_k} } { {\pi x_k} }$ can be taken as $\phi$.
Note that this function $\phi$ is appropriate not only
for the expansions with exact sampled values, but also for falsified expansions.
Indeed, if $N=2$, then $Lf=f+\sum_{[\beta]=1}a_\beta D^\beta f$, where $a_\beta =0$ by~(\ref{fOperL}).
Thus $\phi$ satisfies all conditions of
Theorem~\ref{theoQjAverage2} with  $N=1$,
and, according to~(\ref{fTheoQj32_new}),  the approximation
order of the corresponding  falsified sampling expansions  is $2$
for smooth enough functions $f$.

 Hence, theoretically, we have a simple and very good solution to our problem.
However such expansions are not good from the computational point of view because
in the case $Lf\ne f$ we will not be able to derive explicit formulas
 which is needed for implementations. $\Diamond$

\medskip

II. Let
	$$
\h\phi(\xi) =\prod_{k=1}^d\left(\frac {\sin  {\pi \xi_k} } { {\pi \xi_k} }\right)^2.
$$
Since $\phi(x)=\prod_{k=1}^d(1-|x_k|)\chi_{[-1,1]^d}(x),$
$\phi$ is compactly supported and in ${\cal L}_{p}$.
Also, the function $\sum_{k\in\zd}  |\h\phi(\xi+k)|^q$  is bounded,
$\h\phi$ is continuously differentiable up to any order, the function
$\sum_{l\in\zd, l\neq 0}  |D^{\beta}\h\phi(\xi+l)|$  is bounded
near  the origin for $[\beta]=2.$
Also, the Strang-Fix condition of order~$2$ holds for $\phi$.
The values of $\h\phi$ and its derivatives at the origin are
	$$
	\h\phi(\nul)=1,\quad D^{\beta}\h\phi(\nul)=0,\quad
	[\beta]=1.
	$$
 So, if $\w\phi=\delta,$  then all assumptions
 of Theorem~\ref{theoQjL} are satisfied.
The corresponding sampling expansion of a signal $f$ interpolates
$f$ at the points $M^{-j}k$, $k\in\zd$, the approximation order
depends on how smooth is  $f$, but, according to~(\ref{fTheoQjL}),
 it cannot be better than $2.$
Again $\w\phi$ is associated with the differential
operator $Lf=f+\sum_{[\beta]=1}a_\beta D^\beta f$, where $a_\beta =0$.
Hence the functions $\phi, \w\phi$ satisfy all conditions of
Theorem~\ref{theoQjAverage} with $n=2$, $N=1$,
and, according to~(\ref{fTheoQj3}), the approximation
order of the corresponding  falsified sampling expansions  is $2$
for smooth enough functions $f$.	$\Diamond$

\medskip

III. Let $d=2$,
	$$\h\phi(\xi_1,\xi_2) =
	\frac {1} {(\pi^2\xi_1 \xi_2)^3}
	\left(\sin^3 \pi \xi_1
	\sin^3 \pi \xi_2 +
	b_1\sin^3 \pi \xi_1
	\sin^4 \pi \xi_2+
	b_2\sin^4 \pi \xi_1
	\sin^3 \pi \xi_2\right).$$
Again $\phi$ is compactly supported and belongs to ${\cal L}_p$.
Also, the function $\sum_{k\in\zd}  |\h\phi(\xi+k)|^q$  is bounded,
$\h\phi$ is continuously differentiable up to any order.
Since the trigonometric polynomial in the numerator of $\h\phi$ is bounded, the function
$\sum\limits_{l\in\zd, l\neq 0}  |D^{\beta}\h\phi(\xi+l)|$  is bounded
near  the origin for $[\beta]=3.$
Also, the Strang-Fix condition of order $3$ holds for $\phi$.
The values of $\h\phi$ and its derivatives at the origin are
	$$
	\h\phi(0,0)=1,\quad
	 D^{(1,0)}\h\phi(0,0)=D^{(0,1)}\h\phi(0,0)=0,
	$$
	$$
	 D^{(2,0)}\h\phi(0,0)=\pi^2(2b_1-1),\quad  D^{(0,2)}\h\phi(0,0)=\pi^2(2b_2-1),\quad
	 D^{(1,1)}\h\phi(0,0)=0.
	$$
Now, we choose an appropriate differential operator $L$ in the form $Lf=f+a_{(2,0)} D^{(2,0)}f+a_{(0,2)} D^{(0,2)}f$, or equivalently,
the associated distribution $\w\phi= \delta + \overline{a_{(2,0)}}D^{(2,0)}\delta+
\overline{a_{(0,2)}} D^{(0,2)}\delta$.
Since
	$\h{\w\phi} (\xi)= 1 - 4\pi^2  \overline{a_{(2,0)}}\xi_1
	-4\pi^2  \overline{a_{(0,2)}}\xi_2,$
we have
	$$
\h{\w\phi} (0,0)=1,\quad D^{(2,0)}\h{\w\phi} (0,0) = -4\pi^2 \overline{a_{(2,0)}},	
	\quad D^{(0,2)}\h{\w\phi} (0,0) =  -4\pi^2  \overline{a_{(0,2)}},
$$
	$$
	 D^{(1,0)}\h{\w\phi}(0,0)=D^{(0,1)}\h{\w\phi}(0,0)= D^{(1,1)}\h{\w\phi}(0,0)=0.
	$$
To satisfy condition $D^{\beta}(1-\h\phi \h{\w\phi})(0,0)=0$ for $[\beta]<3$, we have to provide
	$$
	b_1 =\frac 12 ( 1- 4  \overline{a_{(2,0)}}), \quad
	b_2 = \frac 12 (1 - 4 \overline{a_{(0,2)}}).
	$$
Finally, all conditions of Theorem~\ref{theoQjL} are satisfied.
 The  approximation order depends on how smooth is
 $f$, but, according to~(\ref{fTheoQjL}), it cannot be better than $n=3.$

Now we show that the coefficients $b_1, b_2$ can be chosen
such that all conditions of Theorem~\ref{theoQjAverage} are satisfied
with $n=3$, $N=2$.
In this case the differential operator $L$ is given by~(\ref{fOperL}). The
 coefficients $a_{\beta},$ $[\beta]<3$, are as follows
	$$ a_{0,0}=1, \quad a_{1,0} = a_{0,1}= a_{1,1}= 0, \quad
	a_{2,0}= a_{0,2} = \frac 18 \int\limits_0^\infty h(x)^2 \omega(x) dx.$$
	Thus, we set
	  $b_1 = b_2 = \frac 12 (1-4 a_{(2,0)}).$
According to~(\ref{fTheoQj3}), the approximation
order of the corresponding  falsified sampling expansions  is $3$
for smooth enough functions $f$. $\Diamond$	

\medskip

IV. Let $d=1,$
	$$\h\phi(\xi)=
	\frac {\sin^4 \pi\xi + b_1 \sin^5 \pi \xi +
	b_2 \sin^6 \pi\xi + b_3 \sin^7 \pi\xi}
	{(\pi \xi)^4}.$$
Since $\phi$ is bounded and compactly supported,
it is in ${\cal L}_p$,
and $\h\phi$ is  continuously differentiable up to any order.
Also, $\sum_{k\in\zd}  |\h\phi(\xi+k)|^q$  is bounded.
Since the trigonometric polynomial in the numerator of $\h\phi$ is bounded, the
function $\sum_{l\in\zd, l\neq 0}  |D^{\beta}\h\phi(\xi+l)|$  is bounded
near  the origin for $\beta=4.$
Also the Strang-Fix condition of order $4$ holds for $\phi$.
The values of $\h\phi$ and its derivatives at the origin are
	$$
	\h\phi(0)=1,\quad \h\phi'(0)=b_1 \pi,\quad
	\h\phi''(0)=\frac 23\pi^2 (3 b_2 - 2),\quad
	\h\phi'''(0)=\pi^3 (6 b_3 - 5 b_1).
	$$
Now, we choose the appropriate differential operator $L$ in the form $Lf=f+a_1 f'+a_2 f'' + a_3 f'''$, or equivalently,
the associated distribution $\w\phi= \delta - \overline{a_1}\delta'+
\overline{a_2} \delta'' - \overline{a_3} \delta'''$.
Since
	$\h{\w\phi} (\xi)= 1 - 2\pi i \overline{a_1}\xi
	-4\pi^2 \overline{a_2} \xi^2 + 8 \pi^3 i
	\overline{a_3} \xi^3 ,$
we have
	$$\h{\w\phi} (0)=1,\quad \h{\w\phi}' (0) = -2\pi i  \overline{a_1},	
	\quad \h{\w\phi}'' (0) = -8\pi^2 \overline{a_2},\quad
	\h{\w\phi}''' (0) = 48 \pi^3 i 	\overline{a_3}.$$
To satisfy condition $D^{\beta}(1-\h\phi \h{\w\phi})(0)=0$ for $\beta=0,1,2,3$,
we have to provide
\ba
(1-\h\phi\h{\w\phi})'(0)=\pi (b_1-2 i \overline{a_1})=0,
\nonumber
\\
(1-\h\phi\h{\w\phi})''(0)=\frac 23 \pi^2 (-2+3 b_2  -12 \overline{a_2} - 6 i \overline{a_1} b_1)=0,
\nonumber
\\
(1-\h\phi\h{\w\phi})'''(0)=-\pi^3 (5 b_1 - 6 b_3 + 4 i (3b_2 - 2 )\overline{a_1} + 24 b_1 \overline{a_2} - 48 i \overline{a_3})=0.
\label{fExamCoef}
\ea
Thus, the coefficients of the function $\h\phi$ can be easily  found  using the coefficients of the differential operator $L$.
Finally, all conditions of Theorem~\ref{theoQjL} are satisfied.
 The  approximation order depends on how smooth is
 $f$, but, according to~(\ref{fTheoQjL}) with $|\lambda|=|M|$, it cannot be better than $n=4.$

Now we show that the coefficients $b_1, b_2, b_3$ can be chosen
such that all conditions of Theorem~\ref{theoQjAverage} are satisfied with
$n=4$, $N=3$.
In this case the differential operator $Lf=a_0f+a_1 f'+a_2 f'' + a_3 f'''$ is given by~(\ref{fOperL}) and its
 coefficients are defined as
	$$ a_0=1, \quad a_1 = 0, \quad a_2=  \frac 16 \int_0^\infty h(x)^2 \omega(x) dx, \quad a_3=0.$$
	Using~(\ref{fExamCoef}),
	 we set $b_1 = 0, b_2 =
	\frac 23 + 4 a_2, b_3 = 0.$
According to~(\ref{fTheoQj3}), the approximation
order of the corresponding  falsified sampling expansions  is $4$
for smooth enough functions $f$.

Note that all conditions of Theorem~\ref{theoQjAverage1} are also satisfied,
which provides approximation order for a wider class of functions $f$.
Namely,  according to~(\ref{fTheoQj31}),
the approximation order is $3+\frac1p$, whenever $f\in W_1^4$, $f^{(IV)}\in {\rm Lip}_{L_1}\,\varepsilon$, $\varepsilon>0$.  $\Diamond$


\begin{thebibliography}{99}



\bibitem{Boor}
	{\sc C. de Boor}, Quasiinterpolants and approximation power of multivariate splines, Computation of Curves and Surfaces,
 NATO ASI Series \textbf{307} (1990), 313--345.


\bibitem{Borel}
        {\sc E. Borel},
        Sur l'interpolation,
        {\it C. R. Acad. Sci. Paris} {\bf 124} (1897), 673--676.

\bibitem{Butz4}
        {\sc C. Bardaro, P.\,L. Butzer, R.\,L. Stens, and G. Vinti},
        Approximation error of the Whittaker cardinal series in terms of an
        averaged modulus of smoothness covering discontinuous signals,
        {\it Math. Anal. Appl.} {\bf 316} (2006), no.~1, 269--306.

\bibitem{Butz6}
        {\sc C. Bardaro, P.\,L. Butzer, R.\,L. Stens, and G. Vinti},
        Prediction by samples from the past with error estimates covering
        discontinuous signals,
        {\it IEEE Trans. Inform. Theory} {\bf 56} (2010), no.~1, 614--633.



\bibitem{Brown}
        {\sc J.\,L. Brown, Jr.},
        On the error in reconstructing a non-bandlimited function by means of
        the bandpass sampling theorem,
        {\it J. Math. Anal. Appl.} {\bf 18} (1967), 75--84.

\bibitem{Butz5}
        {\sc P.\,L. Butzer, J.\,R. Higgins, and R.\,L. Stens},
        Classical and approximate sampling theorems: studies in the
        $L_p(\mathbb{R})$ and the uniform norm,
        {\it J. Approx. Theory} {\bf 137} (2005), no. 2, 250--263.

\bibitem{Butz8}
        {\sc P.\,L. Butzer and J. Lei},
Approximation of signals using measured sampled values and error analysis,
\textit{Commun. Appl. Anal.} \textbf{4} (2000), no. 2, 245--255.


 \bibitem{ButzNEW}
 {\sc P. L. Butzer, S. Ries, R.L. Stens}, Approximation of continuous and discontinuous functions by generalized sampling series, {\it J. Approx. Theory} \textbf{50}  (1987), no.~1, 25--39.

\bibitem{04}
        {\sc A.\,L. Cauchy},
        Me'moire sur diverses formules d'analyse,
        {\it C. R. Acad. Sci.} {\bf 12} (1841), 283--298.

\bibitem{ChD}
	{\sc Ch.\,Chui and H.\,Diamond}, A characterization of multivariate quasi-interpolation formulas and its applications, {\it Num. Math.} \textbf{57} (1990), 105--121.


 \bibitem{HanSymSmooth}
{\sc B. Han}, Computing the smoothness exponent of a symmetric multivariate refinable function, \textit{SIAM. J. Matrix Anal. \& Appl.} \textbf{24} (2003), no.~3, 693--714.


\bibitem{JZ}
{\sc K.~Jetter and D.~X.~Zhou}, Order of linear approximation from shift invariant
spaces, \textit{Constr. Approx.} \textbf{11} (1995), no.~4, 423--438.

\bibitem{JZ1}
{\sc K. Jetter and D. X. Zhou},
Order of linear approximation on finitely generated
shift invariant spaces", preprint, (1998).


\bibitem{JiaMicPrewav}
 {\sc R.-Q. Jia, C.A. Micchelli}, Using the refinement equations
for the construction of pre-wavelets II: Powers of two,
P.J. Laurent, A. Le M\'ehaut\'e, L.L. Schumaker (Eds.),
Curves and Surfaces, Academic Press, New York (1991), 209--246.

\bibitem{JiaIsotropic}
 {\sc R.-Q. Jia}, Approximation properties of multivariate wavelets,
 	\textit{Math. Comp.} \textbf{67} (1998), no.~222, 647--665.

 \bibitem{v58}
{\sc R.-Q. Jia},  Refinable shift-invariant spaces: From splines to wavelets.
[CA] Chui, C. K. (ed.) et al., Approximation theory VIII. Vol. 2.
Wavelets and multilevel approximation. Papers from the 8th Texas international
conference, College Station, TX, USA, January 8--12, 1995. Singapore:
World Scientific. Ser. Approx. Decompos. 1995, Vol. 6, 179--208.



\bibitem{Jia2}
{\sc R.-Q. Jia}, Approximation by quasi-projection operators in Besov spaces.
{\it J. Approx. Theory} \textbf{162} (2010), no.~1, 186--200.



\bibitem{02}
        {\sc V.\,A. Kotelnikov},
        On the transmission capacity of the 'ether'and of cables in electrical
        communications,
        {\it in} ``Proc. of the First All-Union Conference on the Technological
        Reconstruction of the Communications Sector and Low-current Engineering'',
        Izd. Red. Upr. Svyazzi RKKA, 1933 (in Russian).

\bibitem{KS}
{\sc A. Krivoshein, M. Skopina}, Approximation by frame-like wavelet systems,
{\it Appl. Comput. Harmon. Anal.}
\textbf{31} (2011), no.~3, 410--428.


\bibitem{KS1}
{\sc A. Krivoshein, M. Skopina}, Multivariate sampling-type approximation,
{\it Anal. Appl.} (2016),
DOI: 10.1142/S0219530516500147.

\bibitem {LJC}
{\sc J. J. Lei, R. Q. Jia, E. W. Cheney}, Approximation from
shift-invariant spaces by integral operators, {\it SIAM J. Math. Anal.}
\textbf{28} (1997), 481--498.

\bibitem{NPS}
 {\sc    I. Ya. Novikov, V. Yu. Protasov, and M. A. Skopina}, Wavelet Theory, Amer. Math. Soc, 2011.

\bibitem{03}
        {\sc K. Ogura},
        On a certain transcendetal integral function in the theory of interpolation,
        {\it T\^ohoku Math. J.} {\bf 17} (1920), 64--72.

        \bibitem{01}
        {\sc C.\,E. Shannon},
        Communication in the presence of noise,
        {\it Proc. IRE} {\bf 37} (1949), 10--21.

		  \bibitem{Si1}
				{\sc W. Sickel}, Some remarks on trigonometric interpolation on the $n$-torus,
				{\it Z. Anal. Anwend.} \textbf{10} (1991), no.~4, 551--562.
				
				 \bibitem{Si2}
				{\sc W. Sickel}, Spline representations of functions in Besov-Triebel-Lizorkin spaces on $\mathbb{R}^n$,
				{\it Forum Math.} \textbf{2} (1990), no.~5, 451--475.
				
\bibitem{Sk1}
         {\sc M. Skopina},
        Band-limited scaling and wavelet expansions,
        {\it Appl. Comput. Harmon. Anal.} {\bf 36} (2014), 143--157.
\bibitem{Sk2}
         {\sc M. Skopina},
        Approximation by Band-limited Scaling and Wavelet Expansions,
        Proceedings of International Conferences
        CONSTUCTIVE THEORY OF FUNCTIONS, Sozopol 2013, dedicated to
				Blagovest Sendov and to the memory of Vasil Popov, Prof. Marin Drinov Academic Publishing House,
				Sofia, 2014, 235--251.
				

\bibitem{Stens}
        {\sc R.\,L. Stens},
        Sampling with generalized kernels,
        {\it in} ``Sampling Theory in Fourier and Signal Analysis: Advanced Topics''
        (J.\,R. Higgins and R.\,L. Stens, Eds.),
        Clarendon Press, Oxford, 1999.

 \bibitem{Tr4}
        {\sc A.\,Yu. Trynin},
        On the divergence of sinc-approximations everywhere on $(0,\pi)$,
                 \textit{St. Petersburg Math. J.} \textbf{22} (2011), no.~4,
                 683--701.

\bibitem{Vladimirov-1}
 {\sc V.~S. Vladimirov},
Generalized functions in mathematical physics,
MIR, 1979 (translated from Russian).

\bibitem{W}
        {\sc E.\,T. Whittaker},
        On the functions which are represented by the expansion of the
        interpolation-theory,
        {\it Proc. Roy. Soc. Edinburgh Sect. A} {\bf 35} (1915), 181--194.

\end{thebibliography}
\end{document}